\documentclass[nonblindrev]{informs-ijoo}

\OneAndAHalfSpacedXI

\usepackage{natbib}
 \bibpunct[, ]{(}{)}{,}{a}{}{,}%
 \def\bibfont{\small}%
\usepackage{lipsum}
\usepackage{amsfonts}
\usepackage{graphicx}
\usepackage{epstopdf}
\usepackage{braket,amsfonts}
\usepackage{amsmath,amssymb}
\usepackage{stmaryrd}
\usepackage{array}
\usepackage{threeparttable, booktabs, caption, makecell, multirow}
\usepackage[caption=false]{subfig}
\usepackage{algorithmic}
\usepackage{algorithm}

\usepackage{hyperref}

\usepackage{amsopn}

\usepackage[capitalize, nameinlink]{cleveref}[0.19]

\newcommand{\N}{\mathbb{N}}
\newcommand{\Z}{\mathbb{Z}}
\newcommand{\R}{\mathbb{R}}
\newcommand{\tf}{\tilde{f}}

\TheoremsNumberedThrough     
\ECRepeatTheorems

\EquationsNumberedThrough    

\MANUSCRIPTNO{MS-0001-1922.65}

\begin{document}


\RUNAUTHOR{K. Sun and A. Sun }

\RUNTITLE{Algorithms for DC Programs Based on DME Smoothing}

\TITLE{Algorithms for Difference-of-Convex (DC) Programs Based on Difference-of-Moreau-Envelopes Smoothing }

\ARTICLEAUTHORS{%
\AUTHOR{Kaizhao Sun }
\AFF{School of Industrial and Systems Engineering, Georgia Institute of Technology, \EMAIL{ksun46@gatech.edu}}
\AUTHOR{Xu Andy Sun}
\AFF{Operations Research and Statistics, Sloan School of Management, Massachusetts Institute of Technology, \EMAIL{sunx@mit.edu}}
} 

\ABSTRACT{%
In this paper we consider minimization of a difference-of-convex (DC) function with and without linear equality constraints. We first study a smooth approximation of {a generic} DC function, {termed difference-of-Moreau-envelopes (DME) smoothing}, where \emph{both} components of {the} DC function are replaced by their respective Moreau envelopes. The resulting smooth approximation is {shown to be Lipschitz differentiable, capture stationary points, local, and global minima of the original DC function, and enjoy some growth conditions, such as level-boundedness and coercivity, for broad classes of DC functions. For a smoothed DC program without linear constraints, it is shown that the classic gradient descent method as well as an inexact variant converge to a stationary solution of the original DC function in the limit with a rate of $\mathcal{O}(K^{-1/2})$}, where $K$ is the number of proximal evaluations of both components. Furthermore, when the DC program is explicitly constrained in an affine subspace, we combine the smoothing technique with the augmented Lagrangian function and derive two variants of the augmented Lagrangian method (ALM), named LCDC-ALM {and composite LCDC-ALM}, targeting on different structures of the DC {objective} function. We show that {both algorithms} find an $\epsilon$-approximate stationary solution {of the original DC program} in $\mathcal{O}(\epsilon^{-2})$ iterations. Comparing to existing {methods} designed for linearly constrained weakly convex minimization, the proposed {ALM-based algorithms} can be applied to a broader class of problems, where the objective contains a nonsmooth concave component. Finally, numerical experiments are presented to demonstrate the performance of the proposed algorithms.
}%



\maketitle

%

\section{Introduction}
In this paper, we consider the following {unconstrained} difference-of-convex (DC) program 
\begin{equation}\label{eq: dc}
	  \min_{x\in \R^n}  F(x) := \phi(x) - g(x),
\end{equation}
and the linearly constrained DC (LCDC) program
\begin{equation}\label{eq: lcdc}
	\min_{x \in \R^n} F(x) \quad \mathrm{s.t.} \quad Ax=b,
\end{equation}
where $\phi:\R^n\rightarrow \overline{\R}:=\R\cup\{+\infty\}$ is a weakly-convex function, $g:\R^n\rightarrow \R$ is a {proper and closed} convex function, $A\in \R^{m\times n}$, and $b\in \R^m$. We also work with a more structured setting, where $\phi = f +h$ in a composite form, $h:\R^{n} \rightarrow \overline{\R}$ is a proper, closed, and convex function, and $f: \R^{n} \rightarrow \R$ has a Lipschitz gradient over the effective domain of $h$. DC functions include some important classes of nonconvex functions, such as twice continuously differentiable functions on {any compact convex subset of} $\R^n$ \citep{TuyBook}, {continuous} piecewise-linear functions \citep{melzer1986expressibility}, and multivariate polynomial functions \citep{ahmadi2018dc}. DC programs \eqref{eq: dc}-\eqref{eq: lcdc} appear in various applications such as compressed sensing \citep{yin2015minimization}, high-dimensional statistical regression \citep{cao2018unifying}, and power allocation in digital communication systems \citep{alvarado2014new}, just to name a few. 

In \eqref{eq: dc} and \eqref{eq: lcdc}, we only require the first objective component $\phi$ to be weakly convex, i.e., there exists $m_\phi \geq 0$ such that $\phi + \frac{m_\phi}{2}\|\cdot\|^2$ is a convex function, where $\|\cdot\|$ denotes the Euclidean norm in $\R^n$. Notice that the objective $F$ is still a DC function since $F=\left( \phi + \frac{m_\phi}{2}\|\cdot\|^2\right) -\left(g+ \frac{m_\phi}{2}\|\cdot\|^2\right)$. sNevertheless, we allow $\phi$ to be merely weakly convex, as a convex decomposition of $F$ may not be readily available or necessary. Moreover, we distinguish the DC program \eqref{eq: dc} from the LCDC program \eqref{eq: lcdc} for algorithmic purposes. 
In particular, directly turning \eqref{eq: lcdc} to \eqref{eq: dc} by incorporating linear constraints into the objective using an indicator function may complicate computation, and in situations such as distributed optimization, the data in the linear constraints may not be available to all agents in the optimization problem.

In this paper, we first study a smoothing technique to obtain Lipschitz differentiable DC approximations of a DC function. Then, we combine the smoothing technique with first-order algorithms to solve the DC programs \eqref{eq: dc} and \eqref{eq: lcdc}. It is well-known that the Moreau envelope of a proper closed convex function is a Lipschitz differentiable convex function. Similarly, a nonsmooth weakly convex function can be smoothed by its Moreau envelope. However, directly applying Moreau envelope to the DC function $F=\phi-g$ as a whole might be problematic. On the one hand, the proximal mapping of $F$ can be difficult to compute or even not well-defined; on the other hand, due to the concave component $-g$, the Moreau envelope of $F$ might not be smooth (see Fig. \ref{fig: nonsmooth ME} for an illustration). With these observations, we propose to study a simple smoothing technique, termed difference-of-Moreau-envelope (DME) smoothing, which replaces (weakly) convex components $\phi$ and $g$ separately by their respective Moreau envelopes, i.e., we form
\begin{align}
F_{\mu} := M_{\mu \phi} - M_{\mu g},\label{eq: smooth approx}
\end{align}
where $M_{\mu f}(z):=\min_{x\in \R^n}\{f(x) + \frac{1}{2\mu}\|x-z\|^2\}$ is a Moreau envelope of $f$ with parameter $\mu> 0$. It is easy to see that $F_{\mu}$ is a Lipschitz {differentiable} DC function for all $\mu>0$. Moreover, it has been shown that $F_{\mu}$ preserves the global minimum of $F$ when it exists, i.e. $\min\ F_{\mu}=\min\ F$. Consequently, it is natural to use $F_\mu$ as a surrogate function in search of a minimizers of $F$.

As we finished the manuscript, we discovered that this smoothing technique has been considered as a realization of the Toland duality by \cite{ellaia1984contribution} and applied to smooth DC functions by \cite{hiriart1985generalized,hiriart1991regularize} as a result. In addition, a very recent paper by \cite{themelis2020new} proposed to run gradient descent on $F_\mu$, i.e., through the update
\begin{align}\label{eq: GD}
	z^{k+1} = z^k - \alpha \nabla F_\mu(z^k)
\end{align}
for some proper stepsize $\alpha >0$, which coincides with one of the four algorithms studied in an early version of the present paper. Nevertheless,  (i) some important properties of the smoothed DC function $F_{\mu}$, such as the correspondence between the local minima of $F$ and $F_{\mu}$ and the growth properties of $F_{\mu}$, (ii) an inexact gradient descent method on $F_\mu$, and (iii) the combination of DME technique with the augmented Lagrangian framework for LCDC \eqref{eq: lcdc} are new to our best knowledge. We summarize our contributions in the next subsection.

\subsection{Contributions}
\begin{enumerate}
\item {We carry out a new study on {the {DME} smoothing} of DC functions using two separate Moreau envelopes of the convex components as in  $F_{\mu}$ defined in \eqref{eq: smooth approx}. We show that local minimal solutions of $F_\mu$ can be mapped to those of $F$, complementing existing results on the correspondence between stationary points and global minima of $F_\mu$ and $F$.} Moreover, we identify four {general} conditions on $F$ such that, under each of them, the smoothed DC function $F_{\mu}$ enjoys coercivity or level-boundedness, which is important for the convergence analysis of the proposed algorithms that use $F_{\mu}$ as a surrogate for solving the DC programs.

\item {We propose gradient-based algorithms on $F_\mu$ for the unconstrained DC program \eqref{eq: dc}. As shown in \cite{themelis2020new} and also independently developed in this paper, the classic gradient descent (GD) algorithm applied to  $F_{\mu}$ converges with rate $\mathcal{O}(K^{-1/2})$ where $K$ is the number of proximal evaluations of $\phi$ and $g$, and finds a stationary point of $F$ in the limit.} When $\phi = f+h$, where $f$ has Lispchitz gradient and $h$ is proper, closed, and convex, we propose an inexact GD method for minimizing $F_{\mu}$ that admits a simpler subproblem: $f$ can be replaced by its linearization during the proximal evaluation of $\phi$. The inexact GD achieves the same $\mathcal{O}(K^{-1/2})$ convergence rate, and can be realized as a variant of the proximal DC Algorithm (pDCA) with a novel choice of the proximal center. {A key feature in contrast to DCA-type algorithms is that, we perform a proximal instead of subgradient evaluation on $g$, and proximal mappings of $h$ (or $\phi$ in GD) and $g$ can be implemented in parallel.}
	
\item For the linearly constrained DC program \eqref{eq: lcdc}, we apply the smoothing property of Moreau envelopes to the classic augmented Lagrangian function in augmented Lagrangian method (ALM) and propose two ALM-base algorithms, LCDC-ALM and {composite LCDC-ALM}, for different realizations of $\phi$ (see Table \ref{table: summary} below). In each iteration of LCDC-ALM, the augmented Lagrangian function is smoothed by the DME technique, and then a primal descent step and a dual ascent step are performed. In composite LCDC-ALM, we replace the negative component $-g$ by its linearization, and then smooth a linearized augmented Lagrangian function with its Moreau envelope. For both algorithms, we prove that an $\epsilon$-stationary solution to \eqref{eq: lcdc} can be found in at most $\mathcal{O}(\epsilon^{-2})$ subproblem oracle calls. Moreover, we establish an $\tilde{\mathcal{O}}(\epsilon^{-3})$ first-order complexity estimate for composite LCDC-ALM when each subproblem is solved by a proper first-order algorithm invoking gradient oracles of $f$ and proximal oracles of $h$. LCDC-ALM and composite LCDC-ALM extend ALM-based algorithms considered in \cite{hajinezhad2019perturbed,hong2016decomposing,LiXu2020almv2,melo2020almfullstep,melo2020iteration,zhang2020proximal} by allowing a concave component $-g$ in the objective, though we believe that the aforementioned algorithms have the potential to handle DC objective functions. 	
\end{enumerate}

In addition to the gradient descent scheme \eqref{eq: GD}, whose theoretical properties have already been established by \cite{themelis2020new}, we present three new algorithms in this paper for solving structured DC programs \eqref{eq: dc} and \eqref{eq: lcdc}, summarized as follows.

\begin{table}[!ht]
\caption{Summary of the Proposed Algorithms}\label{table: summary}
\begin{center}
\begin{tabular}{c@{\hskip 3pt}|c@{\hskip 3pt}|c@{\hskip 3pt}|c@{\hskip 3pt}|c@{\hskip 3pt}|c}
\hline
\multirow{2}{*}{Problem} & \multirow{2}{*}{$f$} & \multirow{2}{*}{$h = \phi - f$} & \multirow{2}{*}{$g$} & \multirow{2}{*}{Algorithm} & $\mathcal{O}(\epsilon^{-2})$ Iter. \\
                      &            &                                &     &         &    Complexity        \\
\hline
  \multirow{2}{*}{\eqref{eq: dc}}   &\multirow{2}{*}{Lip. gradient} & \multirow{2}{*}{proper closed convex}       &  \multirow{2}{*}{$\mathcal{P}$}   &  Inexact GD on $F_\mu$  & \multirow{2}{*}{Theorem \ref{thm: igd}}  \\
 & & & & (Algorithm \ref{alg: inexact gd on F_mu})  & \\
\hline
\multirow{4}{*}{\eqref{eq: lcdc}} & \multirow{2}{*}{Lip. gradient}  & \multirow{2}{*}{0}       &  \multirow{2}{*}{$\mathcal{P}$}   &  \small{LCDC-ALM}    & \multirow{2}{*}{Theorem \ref{thm: smooth lcdc-alm}}        \\
& & & & (Algorithm \ref{alg: smooth lcdc-alm}) & \\
\cline{2-6}
              &\multirow{2}{*}{Lip. gradient}  & convex and Lipschitz    &   \multirow{2}{*}{$\mathcal{G}$}   &   \small{{Composite} LCDC-ALM}     &  \multirow{2}{*}{Theorem \ref{thm: ns LCDC-ALM}}      \\
              &                     & over a compact domain       &     &   (Algorithm \ref{alg: ns-lcdc-alm})        &      \\     
\hline  
\multicolumn{6}{@{}p{1.0\textwidth}@{}}{\footnotesize The function $g$ is assumed to be a {finite-valued} convex function in all three algorithms, where either a proximal ($\mathcal{P}$) or subgradient ($\mathcal{G}$) oracle is used in the proposed algorithms.} 
\end{tabular}
\end{center}
\end{table}

\subsection{Notations and Organization}
We denote the set of positive integers up to integer $p$ by $[p]$, the set of nonnegative integers by $\N$, the set of real numbers and nonnegative real numbers by $\R$ and $\R_+$, respectively, the extended real line by $\overline{\R}:=\R\cup \{+\infty\}$, and the $n$-dimensional real Euclidean space by $\R^n$. For $x, y\in \R^n$, the inner product of $x$ and $y$ is denoted by $\langle x, y \rangle$, and the Euclidean norm of $x$ is denoted by $\|x\|$. For $A \in \R^{m\times n}$, $\|A\|$, $\sigma_{\min}(A)$, $\sigma_{\min}^+(A)$, $\lambda_{\max}(A)$, and $\mathrm{Im}(A)$ denote the matrix norm induced by the Euclidean norm, the smallest singular value, the smallest positive singular value, the largest eigenvalue (when $m=n$), and the column space of $A$, respectively. For $X \subseteq \R^n$, $\delta_X$ denotes its $0/\infty$-indicator function, i.e., $\delta_X(x) = 0$ for $x\in X$, and $\delta_X(x) = +\infty$ for $x \notin X$. For a proper function $\phi:\R^n \rightarrow \overline{\R}$, we denote its effective domain by $\mathrm{dom}~\phi := \{x\in \R^n~|~\phi(x) < +\infty\}$. 

The rest of this paper is organized as follows. Section \ref{sec: relatedwork} reviews related work. Section \ref{sec: DC} introduces the DME technique and studies some new properties of the smooth approximation $F_{\mu}$. Section \ref{sec: DC algorithm} shows that applying GD or inexact GD method to the smooth approximation can deliver a stationary solution of \eqref{eq: dc}. Section \ref{sec: LCDC} presents two ALM-based algorithms with convergence analysis for the LCDC problem \eqref{eq: lcdc} under different assumptions of $\phi$. We present numerical experiments in Section \ref{sec: numerical}, and finally give some concluding remarks in Section \ref{sec: conclusion}.

\section{Related Work}\label{sec: relatedwork}
In this section, we review some related work in smoothing, DC algorithms, and ALM.
\subsection{Smooth Approximation of Nonsmooth Functions}
There is a rich literature on approximation of nonsmooth nonconvex functions. \cite{chen1996class} proposed a class of smooth approximations for the max function $(t)_+=\max\{0,t\}$ using integral convolution with some symmetric and piecewise-continuous density function, which is also closely related to mollification in functional analysis. \cite{chen2012smoothing} exploited this scheme in approximating composition of $C^1$ functions with $(t)_+$. \cite{lu2014iterative} proposed a Lipschitz continuous $\epsilon$-approximation of the function $\|\cdot\|_p^p$ for $p\in (0,1)$, where the tolerance $\epsilon>0$ controls the approximation error. The smooth approximation \eqref{eq: smooth approx} is a different smoothing technique, relying on inf-convolution, rather than on integral convolution. In \cite{ellaia1984contribution}, it is shown that $F_{\mu}$ is a special way to realize Toland duality of DC programs by adding ${\frac{1}{2\mu}}\|\cdot\|^2$ to the convex components. In \cite{hiriart1985generalized} and \cite{hiriart1991regularize}, it is shown that, as a result of Toland duality, the stationary points and global optima of $F_{\mu}$ and $F$ are closely related.

\subsection{DC Algorithms}
Algorithms for minimizing a weakly convex function have been studied in the literature \citep{davis2018stochastic,drusvyatskiy2017proximal}.
When $-g$ cannot be convexified through the addition of a quadratic form, a classic iterative approach to finding stationary solutions of \eqref{eq: dc} is the so-called DC algorithm (DCA) proposed by \cite{tao1997convex}, (see also \cite{tao2005dc}), where $-g$ is replaced by a linear over-approximation in each iteration, using subgradient information; see \cite{artacho2018accelerating,le2012exact,le1997solving,tao1998dc} for specific applications and convergence results of DCA. 

An important variant of DCA is the proximal DCA (pDCA), which, in addition to replacing the concave function with a linear majorization, includes a proximal term in the minimization problem at iteration $k+1$:
	\begin{equation}\label{eq: pDCA}
		x^{k+1} = \argmin_{ x\in \R^n} \left\{\phi(x) -\langle \xi_g^k, x-x^k\rangle +\frac{1}{2c_k}\|x-x^k\|^2\right\},
	\end{equation}
	where $c_k>0$ and $\xi_g^k \in \partial g(x^k)$, e.g. see \cite{souza2016global,sun2003proximal}. When $\phi = f+h$,  
	where $f$ has Lipschitz gradient and $h$ is proper closed and convex, one can further linearize $f$ \citep{an2017convergence}, 
	and replace the proximal center $x^k$ by a proper extrapolation \citep{wen2018proximal}. More recently,
	\cite{de2019proximal,de2020sequential} replaced the positive component $\phi$ with some minorant function such as a 
	bundle of cutting planes. 	
	
\cite{pang2017computing} characterized different types of stationary points and their relations for generic DC programs. They considered the case $g(x) = \max_{i\in I}g_i(x)$, where each $g_i$ is {continuously} differentiable and convex and $I$ is a finite index set, and proposed an enhanced DCA that subsequentially converges to a d(irectional)-stationary point. Furthermore, they showed that the algorithm can be extended to compute a B(ouligand)-stationary point for a class of DC constrained DC program under a suitable constraint qualification. This work also motivates several follow-up works. \cite{lu2019enhanced} extended the problem setting to allow an infinite supremum in the definition of $g$; by incorporating the extrapolation technique into the enhanced DCA in \cite{pang2017computing}, they proposed an algorithm named EPDCA and established an iteration complexity of $\mathcal{O}(\epsilon^{-2})$ for computing an approximate stationary solution with tolerance $\epsilon>0$. \cite{lu2019nonmonotone} later applied nonmonotone line-search schemes and randomized update to accelerate the convergence of EPDCA in practice. \cite{taolinear} established linear convergence of the enhanced DCA by utilizing locally linear regularity and error bound conditions.
	
Another related work by \cite{banert2019general} applied the proximal alternating linearized minimization (PALM) algorithm \citep{bolte2014proximal} to a primal-dual formulation of a DC program. Convergence rate results are established under the Kurdyka–\L ojasiewicz (K\L) property \citep{kurdyka1998gradients,lojasiewicz1963propriete}.

\subsection{ALM for Linearly-constrained Weakly Convex Minimization}
The augmented Lagrangian method (ALM), which was proposed in the late 1960s by \cite{hestenes1969multiplier} and \cite{powell1967method}, provides a powerful algorithmic framework for constrained optimization. The convergence, local, and global convergence rate of ALM have been extensively studied for convex programs \citep{AybatIyengar2013alm,lanMonteiro2016alm,LiQu2019alm,LiuLiuMa2019alm,LuZhou2018alm,rockafellar1973multiplier,rockafellar1976augmented,Xu2021alm} and smooth nonlinear programs \citep{bertsekas2014constrained,Sahin2019alm}. In the following, we review some  recent developments in ALM-based algorithms applied to linearly-constrained nonconvex problems of the form
	\begin{equation}\label{eq: lc-alm-literature}
		\min_{x\in \R^n} \phi(x) =  f(x) +h(x) \quad \mathrm{s.t.}\quad Ax=b,
	\end{equation}
    where $f$ has Lipschitz gradient and $h$ is a possibly nonsmooth convex function. 
    
\cite{hong2016decomposing} considered the case where $h=0$, and proposed a proximal ALM: in iteration $k+1$, an additional proximal term $\frac{\rho}{2}\|x-x^k\|^2_{B^\top B}$ is added to the augmented Lagrangian (AL) function.
Assuming $h$ is a compactly supported convex function, \cite{hajinezhad2019perturbed} proposed a perturbed proximal ALM that will converge to a solution with bounded infeasibility; when the initial point is feasible, they established an iteration complexity of $\mathcal{O}(\epsilon^{-4})$, where $\epsilon>0$ measures both stationarity and feasibility. 
Under the same perturbed AL framework, \cite{melo2020iteration} applied an accelerated composite gradient method (ACG) \citep{beck2009fast} to solve each proximal ALM subproblem and obtained an improved iteration complexity of $\mathcal{O}(\log(\epsilon^{-1})\epsilon^{-3})$ total ACG iterations, which can be further reduced to $\mathcal{O}(\log(\epsilon^{-1})\epsilon^{-2.5})$ with mildly stronger assumptions. In \cite{melo2020almfullstep}, this inner acceleration scheme is embedded in the proximal ALM with full dual multiplier update. 
\cite{LiXu2020almv2} combined an inexact ALM and a quadratic penalty method to solve a convex-constrained program with a weakly-convex objective function, and they showed that an $\epsilon$-KKT solution can be found in $\mathcal{O}(\log(\epsilon^{-1})\epsilon^{-2.5})$ adaptive accelerated proximal gradient steps.
 
Finally we highlight two recent works by \cite{zhang2020global,zhang2020proximal} and a further generalization by \cite{zeng2021moreau}. Zhang and Luo studied the case where $h$ is an indicator function of a hypercube or a polyhedron. In practice, their specially chosen proximal term prevents iterates from oscillation, therefore stabilizing the algorithm; in the convergence proof, the authors constructed a novel potential function, which combines properties from primal descent, dual ascent, and proximal descent.  Utilizing an error bound that exploits the polyhedral structure of $h$, the authors showed the proposed algorithm will find an $\epsilon$-approximate KKT point in $\mathcal{O}(\epsilon^{-2})$ iterations. Moreover, their methodology has been successfully applied to min-max problems \citep{zhang2020single} and distributed optimization \citep{chen2021communication}. More recently, \cite{zeng2021moreau} generalized the setting \eqref{eq: lc-alm-literature} to any weakly-convex objective, and proposed a Moreau Envelope ALM (MEAL) that achieves the same iteration complexity under either the {implicit Lipschitz subgradient} condition or the {implicit bounded subgradient} condition. The novel conditions relax the standard smoothness and bounded subgradient conditions commonly adopted in weakly-convex minimization, and provide a new way to control dual iterates using primal iterates in the analysis of ALM.

\section{Smooth Approximation of DC Functions}\label{sec: DC}
\subsection{Assumptions and Stationarity}
We make the following assumptions in this paper regarding the DC function $F = \phi - g$.
\begin{assumption}\label{assumption: dc}
	The function $\phi:\R^n\rightarrow \overline{\R}$ is proper, lower semi-continuous, and $m_\phi$-weakly convex, i.e., there exists a finite $m_{\phi}\geq 0$ such that $\phi +\frac{m_{\phi}}{2}\|\cdot\|^2$ is a convex function. The function $g:\R^n \rightarrow \R$ is convex. Moreover, when $\phi = f+ h$ admits a composite form, the function $h:\R^n \rightarrow \overline{\R}$ is proper, closed, and convex, and $f:\R^n \rightarrow \R$ is Lipschitz differentiable with modulus  $L_f> 0$ over $\mathcal{H}:=\mathrm{dom}~h =  \{x\in \R^n~|~h(x) < +\infty\}$. 
\end{assumption}
\begin{assumption}\label{assumption: F finite optimal}
	The set of global minimizers of $F$, $\argmin F$, is nonempty.
\end{assumption}
By Assumptions \ref{assumption: dc} and \ref{assumption: F finite optimal}, it holds that 
\begin{align}
    \R^n = \mathrm{dom} ~ g = \mathrm{dom}~\partial g: = \{x\in \R^n: \ \partial g(x) \neq \emptyset \}, ~\text{and}~F^* := \min_{x\in \R^n} F(x) < +\infty.
\end{align}
In addition, we adopt the following definition for stationary points of $F$.
\begin{definition}\label{def: dc stationary}
    We say $x\in\R^n$ is a stationary point of $F$ if
    \begin{equation}\label{eq: dc stationary}
        0 \in \partial \phi(x) -  \partial g(x);
    \end{equation}
    or equivalently, $\partial \phi(x) \cap \partial g(x) \neq \emptyset$. Furthermore, given $\epsilon >0$, we say $x\in\R^n$ is an $\epsilon$-stationary point of $F$ if there exists $(\xi,y)\in \R^n \times \R^n$ such that 
    \begin{align}\label{eq: approx dc stationary}
        \xi \in \partial \phi(x) -  \partial g(y), \ \text{and} \ \max \{\|\xi\|, \|x-y\|\} \leq \epsilon.
    \end{align}	
    When $\phi = f + h$, we simply replace $\partial \phi(x)$ by $\nabla f(x) + \partial h(x)$ in \eqref{eq: dc stationary} and \eqref{eq: approx dc stationary}.
\end{definition}
\begin{remark} 
	We give some comments regarding Definition \ref{def: dc stationary}.
	\begin{enumerate}
		\item We use $\partial \phi(x)$ to denote the general subdifferential of $\phi$ at $x$ \cite[Definition 8.3]{rockafellar2009variational}. If $\phi$ is continuously differentiable in a neighborhood of $x$, then $\partial \phi(x) = \{\nabla \phi(x)\}$; if $\phi = f + h$ with $f$ continuously differentiable and $h$ finite at $x$, then $\partial \phi(x) = \nabla f(x) + \partial h(x)$ \cite[Exercise 8.8]{rockafellar2009variational}.
        \item Compared to \eqref{eq: dc stationary}, a more natural stationarity condition of $F$ would be $0 \in \partial (\phi-g) (x)$. In view of the previous comment, if $g$ is continuously differentiable over $\mathrm{dom}~\phi$, then  $\partial (\phi - g) (x) = \partial \phi(x) -  \nabla g(x)$ for all $x\in \mathrm{dom}~\phi$. However, the equality does not hold in general. The condition \eqref{eq: dc stationary} is known as criticity in the DC literature. Some other stationary conditions include d(irectionary)-stationarity and c(larke)-stationarity, which are defined as the directional derivative or Clarke directional derivative being nonnegative in every feasible direction. Under suitable settings, e.g., see \cite{pang2017computing}, local minimum of $F$ implies
			\begin{align*}
			 \partial g({x}) \subseteq \partial \phi({x})\equiv \text{d-stationary}\Rightarrow \text{c-stationary} \Rightarrow \eqref{eq: dc stationary}.
			\end{align*}	
			Clearly the above three types of stationarity are equivalent when $\partial g({x}) = \{\nabla g({x})\}$. However, since we do not impose any smoothness structure on $g$, computing (or even verifying) a d-stationary point can be intractable, so we adopt the more general notion \eqref{eq: dc stationary} in this paper. More details on different stationarity conditions can be found in \cite{pang2017computing} and \cite{de2019proximal}.
		\item Usually $\partial \phi$ and $\partial g$ are not evaluated at the same point during DCA or pDCA. For example in \eqref{eq: pDCA}, a subgradient $\xi_g^k\in \partial g(x^k)$ is used to set $x^{k+1}$. So \eqref{eq: approx dc stationary} is a natural relaxation of \eqref{eq: dc stationary}, and can serve as the stopping criteria for some iterative approach in practice.
	\end{enumerate}
\end{remark}

\subsection{New Properties of the Difference-of-Moreau-Envelopes Smoothing}
\begin{figure*}[h!]
\caption{An Example: $F=\phi-g$, where $\phi(x)=\delta_{[-1,1]}(x)$, and $g(x)=\frac{1}{2}x^2$}\label{fig: nonsmooth ME}   
\center{
	\begin{tabular}{@{}c@{}}
    	\includegraphics[width=.4\linewidth]{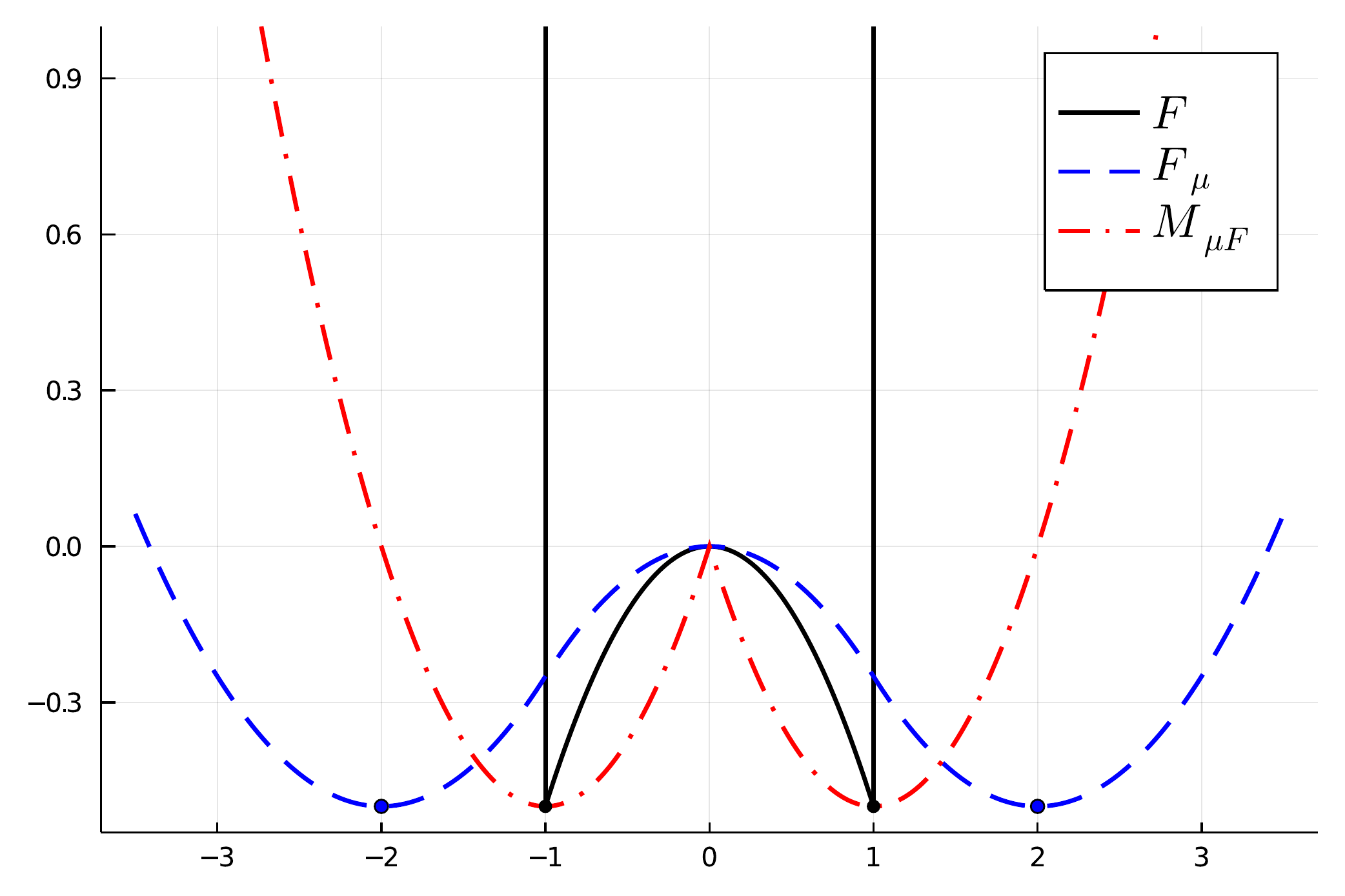} \\[\abovecaptionskip]
    	(a) $F$, $F_\mu$, and $M_{\mu F}$ ($\mu=1$)
  	\end{tabular}
  	\begin{tabular}{@{}c@{}}
    	\includegraphics[width=.4\linewidth]{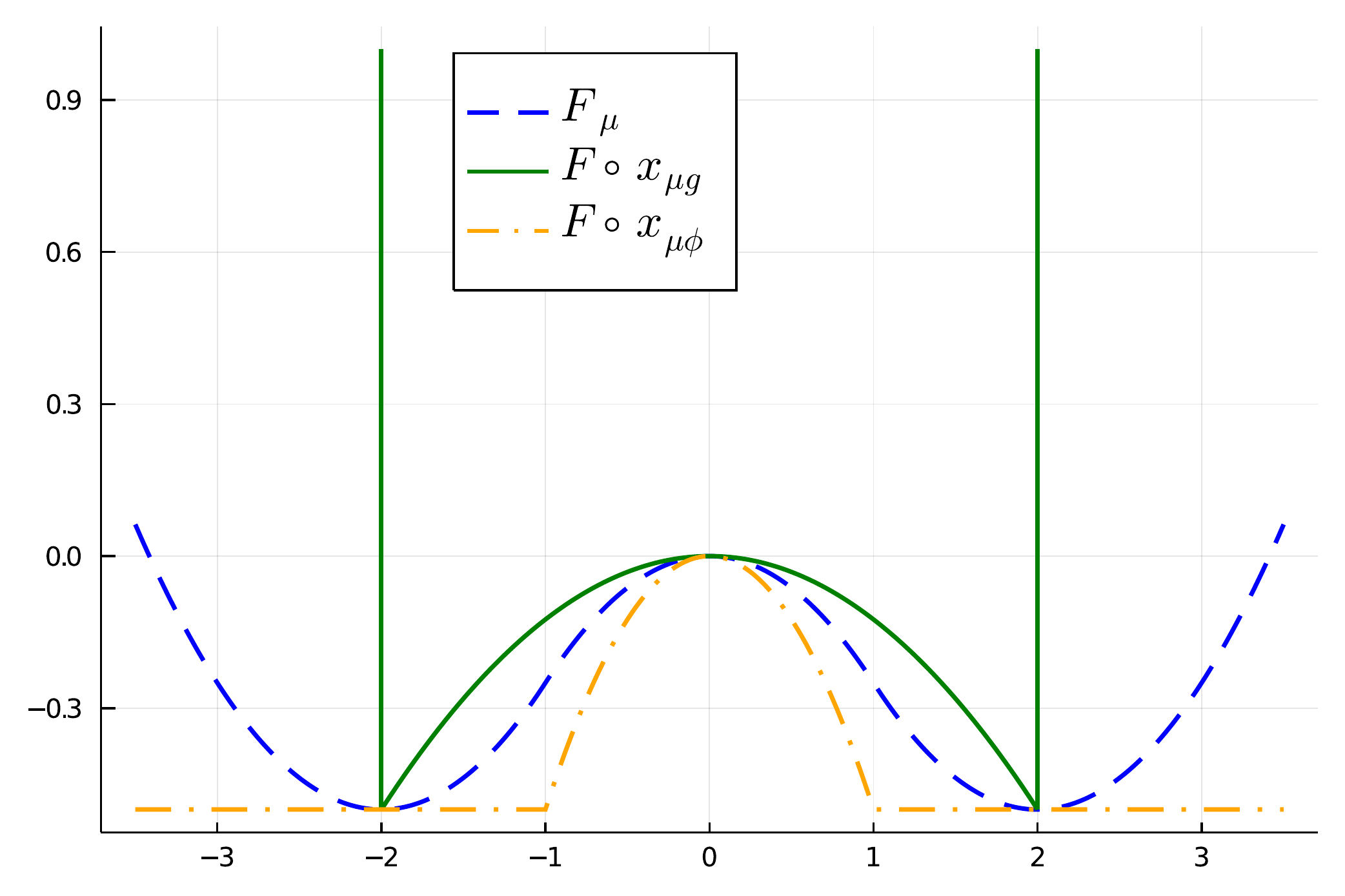} \\[\abovecaptionskip]
    	(b) $F\circ x_{\mu \phi}\le F_\mu\le F\circ x_{\mu g}$ ($\mu=1$)
  	\end{tabular}}
\end{figure*}

Given $0 < \mu < 1/m_{\phi}$, the Moreau envelope $M_{\mu\phi}:\R^n\rightarrow\R$ and the proximal mapping $x_{\mu \phi}:\R^n\rightarrow\R^n$ of $\phi$ are defined as
\begin{align}\label{eq: phi moreau}
	M_{\mu \phi}(z) := \min_{x\in \R^n}\left\{\phi(x) + \frac{1}{2\mu}\|x-z\|^2\right\}, \quad x_{\mu \phi}(z) := \argmin_{x\in \R^n}\left\{\phi(x) + \frac{1}{2\mu}\|x-z\|^2\right\}. 
\end{align}
It is known that $M_{\mu \phi}$ forms a smooth approximation of the possibly nonconvex nonsmooth function $\phi$. Similarly, the Moreau envelope and proximal mapping of $g$ are given by $M_{\mu g}$ and $x_{\mu g}$, respectively. Consequently, the function $F_{\mu}(z) = M_{\mu \phi}(z)-M_{\mu g}(z)$ constitutes a smooth approximation of $F$. Our motivation to form the Moreau envelopes of $\phi$ and $g$ separately comes from the observation that the Moreau envelope of a concave function may not be smooth, therefore, the Moreau envelope {$M_{\mu F}$} of the DC function $F=\phi-g$ as a whole may not be smooth. But smoothing each component of $F$ separately will surely give a smooth DC function. This is shown in the next example. Consider $F=\phi-g$, where $\phi(x)=\delta_{[-1,1]}(x)$, and $g(x)=\frac{1}{2}x^2$. In Fig. \ref{fig: nonsmooth ME}(a) , we plot $F$, the smooth approximation $F_\mu$, and the Moreau envelope $M_{\mu F}$ of $F$; we see that $M_{\mu F}$ is not smooth at the origin. In Fig. \ref{fig: nonsmooth ME}(b), we see the smooth approximation $F_\mu$ is further bounded by $F\circ x_{\mu \phi}$ from below and by $F\circ x_{\mu g}$ from above; this is formally stated in Lemma \ref{lemma: bounds on F_mu} and the proof follows directly from the definitions of $F_\mu$, $x_{\mu g}(z)$, and $x_{\mu \phi}(z)$.

\begin{lemma}[Bounds on $F_\mu$]\label{lemma: bounds on F_mu}    
    Suppose Assumption \ref{assumption: dc} holds and $0 < \mu < 1/m_{\phi}$ in \eqref{eq: phi moreau}. For all $z\in\R^n$, it holds $(F\circ x_{\mu\phi})(z)\le F_\mu(z)\le (F\circ x_{\mu g})(z)$. 
\end{lemma}

In \ref{section: preliminaries}, we provide some preliminary properties on DME smoothing. In particular, it has been established that (i) $F_\mu$ has a Lipschitz gradient with 
\begin{align}\label{eq: gradient of F_mu}
	\nabla F_{\mu}(z) = \mu^{-1}(x_{\mu g}(z)-x_{\mu \phi}(z) )	
\end{align}
for all $z\in \R^n$, (ii) if $z\in \R^n$ is a stationary point or a global minimizer of $F_\mu$, then $x_{\mu \phi}(z) = x_{\mu g}(z)$ is a stationary point or a global minimizer of $F$, respectively, and (iii) the converse of (ii) holds as well. Moreover, when $g$ is smooth or has bounded subgradients, we show below that a local minimum of $F_\mu$ can be used to construct a local minimum of $F$.
\begin{proposition}[Correspondence of Local Minima]\label{prop: map local}
	 Suppose Assumptions \ref{assumption: dc} and \ref{assumption: F finite optimal} hold, $0 < \mu < 1/m_{\phi}$ in \eqref{eq: phi moreau}, and  $\bar{z}$ is a local minimizer of $F_\mu$, i.e., there exists $r>0$ such that $F_\mu(\bar{z}) \leq F_{\mu}(z)$ for all $z\in \R^n$ satisfying $\|z-\bar{z}\|\leq r$. The following two claims hold.
	 \begin{enumerate}
	 	\item Suppose $g$ is Lipschitz differentiable with modulus $L_g$, then $F(x_{\mu \phi}(\bar{z}))\leq F(x)$ for all $x\in \R^n$ such that $\|x-x_{\mu \phi}(\bar{z})\|\leq \frac{r}{\mu L_g+1}$.
	 	\item Suppose $\sup_{\xi_g\in \partial g(\R^n)}\|\xi_g\|\leq M_{\partial g}$ for some $M_{\partial g}\geq0$ and $r > 2\mu M_{\partial g}$, then $F(x_{\mu \phi}(\bar{z}))\leq F(x)$ for all $x\in \R^n$ such that $\|x-x_{\mu \phi}(\bar{z})\|\leq r-2\mu M_{\partial g}.$
	 \end{enumerate}
\end{proposition}
\proof{Proof.}
	See \ref{proof: map local}
\endproof

As a consequence of Proposition \ref{prop: known} and Proposition \ref{prop: map local}, the smooth approximation $F_{\mu}$ can serve as a surrogate for the nonsmooth DC function $F$ when solving problem \eqref{eq: dc}. Before making this idea more concrete in Section \ref{sec: DC algorithm}, we consider an important class of properties defined as follows.
\begin{definition}[Theorem 3.26 of \cite{rockafellar2009variational}]\label{def: level bounded}
	A function $F:\R^n\rightarrow \overline{\R}$ is said to be (i) \textit{level-bounded} if its lower level set $\mathrm{lev}_{\alpha} F:= \{x:F(x)\leq \alpha\}$ is bounded (possibly empty) for any $\alpha \in \R$, (ii) \textit{level-coercive} if there exist some $\alpha \in (0, +\infty)$ and $r \in \R$ such that $F\geq \alpha \|\cdot\|+r$, and (iii) \textit{coercive} if for any $\alpha \in (0, +\infty)$, there exists $r \in \R$ such that $F \geq \alpha\|\cdot\|+r$.
\end{definition}

Clearly, coercivity implies level-coercivity, which in turn implies level-boundedness. Level-boundedness is important for a descent algorithm to produce bounded iterates, so that at least one limit point exists and we can study the asymptotic behavior of the algorithm. We note that for a general DC function $F$ satisfying Assumption \ref{assumption: dc}, $F_\mu$ is not necessarily level-bounded. In the next proposition, we give some sufficient conditions that ensure $F_\mu$ is coercive or level-bounded.
\begin{proposition}\label{prop: level boundedness}
	Suppose Assumptions \ref{assumption: dc} and \ref{assumption: F finite optimal} hold, and $0 < \mu < 1/m_{\phi}$.
	\begin{enumerate}
		\item If $F$ is level-bounded, and $g$ is Lipschitz up to some constant, i.e., there exist $L>0$ and $M\geq 0$ such that $\|g(x)-g(z)\|\leq L\|x-z\|+M$ for all $x,z\in \R^n$, then $F_\mu$ is level-bounded.
 		\item If $F$ is level-coercive, then $F_\mu$ is level-bounded.
		\item If $\mathrm{dom} ~\phi := \{x: \phi(x) < +\infty\}$ is compact, then $F_\mu$ is coercive, and therefore level-bounded.
		\item If $F$ is level-bounded, and $\mathrm{dom}~\phi =\R^n$, then $F_{\mu}$ is level-bounded.
	\end{enumerate}
\end{proposition}
\proof{Proof.}
	See \ref{proof: level boundedness}. \Halmos
\endproof

\section{Unconstrained DC Optimization}\label{sec: DC algorithm}
In this section, we present two algorithms for solving the unconstrained DC program \eqref{eq: dc}. Although we call problem \eqref{eq: dc} ``unconstrained", under Assumption \ref{assumption: dc}, implicit convex functional constraints can be incorporated into the definition of $\phi$ using an indicator function, as long as proximal evaluation of the resulting $\phi$ is available. 

Minimizing $F$ is in general challenging due to its nonconvexity and nonsmoothness. In contrast, the {DME} smooth approximation $F_{\mu}$ provides an attractive surrogate: the Lipschitzness of $\nabla f$ is a desirable property for a wide range of first-order methods, among which the GD \eqref{eq: GD} is probably the most classic one. By \eqref{eq: gradient of F_mu}, each gradient update \eqref{eq: GD} requires the proximal evaluations of $\phi$ and $g$, which can be performed in parallel and is hence a major difference from DCA and pDCA. The convergence of GD \eqref{eq: GD} is stated in the next theorem. A similar result has been proved in \cite[Theorem 7]{themelis2020new}. We provide a detailed proof for self-consistency.

\begin{theorem}\label{thm: gd}
	Suppose Assumptions \ref{assumption: dc} and \ref{assumption: F finite optimal} hold and let $\epsilon>0$, $0 < \mu < 1/m_{\phi}$, and $0 < \alpha \leq  1/L_{F_\mu}$ where $L_{F_\mu} = \frac{2-\mu m_{\phi}}{\mu-\mu^2m_{\phi}}$ if $m_{\phi} >0$ and $L_{F_\mu} =2\mu^{-1}$ if $m_{\phi}= 0$. Let $\{(z^k$, $x_{\mu \phi}(z^k)$, $x_{\mu g}(z^k))\}_{k\in \N}$ be the sequence generated by \eqref{eq: GD} Denote $\xi^k = \mu^{-1}(x_{\mu g}(z^k)$ $- x_{\mu \phi}(z^k)).$ Then for any positive integer $K$, there exists $0\leq \bar{k}\leq K-1$ such that 
	\begin{subequations} \label{eq: gd_opt_coniditon}
		\begin{align}
			& \xi^{\bar{k}} \in \partial \phi(x_{\mu \phi}(z^{\bar{k}})) - \partial g(x_{\mu g}(z^{\bar{k}})), \ \text{and} \\
			& \max\{\|\xi^{\bar{k}}\|, \|x_{\mu g}(z^{\bar{k}}) - x_{\mu \phi}(z^{\bar{k}})\| \}\leq  \max\{1, \mu^{-1}\}  \left(\frac{2\mu^2 (F_{\mu}(z^0)-F^*)}{\alpha K}\right)^{1/2}.
		\end{align}
	\end{subequations}
	Therefore, the GD \eqref{eq: GD} finds an $\epsilon$-stationary solution in the sense of \eqref{eq: approx dc stationary} in no more than 
	\begin{align}\label{eq: gd iter bound}
		  \left \lceil  \frac{ 2\mu^2\max\{1, \mu^{-1}\}^2(F_{\mu}(z^0)- F^*)}{\alpha \epsilon^2} \right \rceil =\mathcal{O}(\epsilon^{-2})
	\end{align}
	iterations. Moreover, if $F_\mu$ is level-bounded, then the sequence $\{z^k\}_{k\in \Z_+}$ stays bounded; for every limit point $z^*$, it holds $x_{\mu \phi}(z^*) = x_{\mu g}(z^*)$, and $x_{\mu \phi}(z^*)$ satisfies the exact stationary condition \eqref{eq: dc stationary}.
\end{theorem}
\proof{Proof.}
	See \ref{appendix: proof_of_gd}. \Halmos
\endproof
Next we consider a more structured setting where $\phi = f+ h$ admits a composite form. Instead of directly evaluating the proximal mapping of $\phi$, which might still be difficult, we can replace the smooth component $f$ with its linearization at the previous iterate. See Algorithm \ref{alg: inexact gd on F_mu}.
\begin{algorithm}[tbh]
	\caption{: Inexact GD on $F_\mu$} \label{alg: inexact gd on F_mu}
	\begin{algorithmic}[1]
		\STATE \textbf{Let} $0 < \mu < 1/L_f$ and $0 < \beta  < 2$;
		\STATE \textbf{Initialize} $x^0, z^0\in \R^n$;
		\FOR{$k = 0,1,\cdots$}
		\STATE $\displaystyle{x^{k+1} = \argmin_{x\in\R^n} \left\{\langle \nabla f(x^k), x \rangle + h(x) + \frac{1}{2\mu}\|x-z^k\|^2\right\}}$;\label{eq: linear prox}
		\STATE $\displaystyle{x_{\mu g}(z^{k}) = \argmin_{x\in\R^n}\left\{g(x) + \frac{1}{2\mu}\|x-z^k\|^2\right\}}$;
		\STATE $z^{k+1} = z^k + \beta (x^{k+1}-x_{\mu g}(z^k))$; \label{eq: inexact gd}
		\ENDFOR
	\end{algorithmic}
\end{algorithm}

\begin{remark}~
	\begin{enumerate}
		\item The update of $x^{k+1}$ in line \ref{eq: linear prox} takes the form 
		$x^{k+1} = \argmin_{x} h(x) + \frac{1}{2\mu}\|x-(z^k-\mu \nabla f(x^k))\|^2.$
		We first take a gradient (forward) step with respect to the smooth component $f$, and then take a proximal (backward) step with respect to the nonsmooth component $h$. A notable difference is that, we move $z^k$, instead of $x^k$, along the direction $-\nabla f(x^k)$.
		\item Our choice of the proximal center $z^k$ is different from the extrapolation $x^{k}+\beta(x^k-x^{k-1})$ considered in \cite{wen2018proximal}. 
		\item We view $\mu^{-1}(x^{k+1}-x_{\mu g}(z^k))$ as an inexact gradient of $F_{\mu}$ at $z^k$, and take a similar gradient step to update $z^{k+1}$ with some properly chosen step size $\beta \mu$ as in line \ref{eq: inexact gd}. 
	\end{enumerate}
\end{remark}

We propose the following potential function
\begin{align}
	\mathcal{F}(x, z) := & \phi(x) + \frac{1}{2\mu}\|x-z\|^2 - M_{\mu g}(z)= f(x)+h(x) + \frac{1}{2\mu}\|x-z\|^2 - M_{\mu g}(z). 
\end{align}
We show the sequence $\{\mathcal{F}(x^k,z^k)\}_{k\in \N}$ generated by Algorithm \ref{alg: inexact gd on F_mu} is bounded and non-increasing.
\begin{lemma}\label{lemma: igd descent}
	Suppose Assumptions \ref{assumption: dc} and \ref{assumption: F finite optimal} hold. For all $k\in \N$, $\mathcal{F}(x^{k}, z^k) \geq F_{\mu}(z^k) \geq F^*$, and 
	\begin{equation}\label{eq: igd descent}
		\mathcal{F}(x^{k}, z^k) -\mathcal{F}(x^{k+1},z^{k+1}) \geq c_1\|x^{k+1}-x^{k}\|^2 + c_2\|z^{k+1}-z^k\|^2,
	\end{equation} 
	where 
	\begin{align}\label{eq: c_12}
		c_1 =  \frac{\mu^{-1}-L_f}{2} >0 , \quad  c_2 =  \frac{1}{\mu}\left(\frac{1}{\beta}-\frac{1}{2}\right)>0.
	\end{align}
\end{lemma}
\proof{Proof.}
	See \ref{proof: igd descent}. \Halmos
\endproof

\begin{theorem}\label{thm: igd}
	Suppose Assumptions \ref{assumption: dc} and \ref{assumption: F finite optimal} hold, and let $\epsilon>0$. Let $\{(x^{k+1}$, $x_{\mu g}(z^k)$, $z^{k+1})\}_{k\in \N}$ be the sequence generated by Algorithm \ref{alg: inexact gd on F_mu}. Denote $\xi^{k+1} := \nabla f(x^{k+1})-\nabla f(x^k) -{\mu}^{-1}(x^{k+1}-x_{\mu g}(z^k))$. Then for any positive integer $K$, there exists $0\leq \bar{k}\leq K-1$ such that
	\begin{subequations}
	\begin{align}
		& \xi^{\bar{k}+1} \in \partial \phi(x^{\bar{k}+1}) - \partial g(x_{\mu g}(z^{\bar{k}})), \ \text{and} \label{eq: igd xi}\\
		& \max\{\|\xi^{\bar{k}+1}\|, \|x_{\mu g}(z^{\bar{k}}) - x^{k+1}\| \}\leq  \left(L_f + \frac{\mu^{-1}+1}{\beta}\right) \left( \frac{\mathcal{F}(x^0,z^0)-F^*)}{\min \{c_1,c_2\}K}\right)^{1/2}.\label{eq: igd max}
	\end{align}
	\end{subequations}
	Therefore, Algorithm \ref{alg: inexact gd on F_mu} finds an $\epsilon$-stationary solution in the sense of \eqref{eq: approx dc stationary} in no more than 
	\begin{align}\label{eq: igd iter bound}
		\left \lceil  \frac{ (L_f+(\mu^{-1}+1)/\beta)^2(\mathcal{F}(x^0, z^0)- F^*)}{\min \{c_1, c_2\}\epsilon^2} \right \rceil =\mathcal{O}(\epsilon^{-2})
	\end{align}
	iterations. Moreover, if $F_\mu$ is level-bounded, then the sequence $\{(x^k,z^k)\}_{k\in \Z_+}$ stays bounded; for every limit point $(x^*, z^*)$, it holds $x^* = x_{\mu g}(z^*)$, and $x^*$ satisfies the exact stationary condition \eqref{eq: dc stationary}.
\end{theorem}	
\proof{Proof.}
	See \ref{proof: igd}.  \Halmos
\endproof
{We end this section with some additional remarks regarding comparisons with (proximal) DCA-type algorithms. Admittedly, (proximal) DCA is also able to achieve a similar iteration complexity, and usually a subgradient evaluation of $g$ is easier than its proximal evaluation.
We note that the DME technique provides a globally smooth approximation for any generic DC function with desirable properties established in Section \ref{sec: DC}, which a local linearization of $g$ used in DCA does not possess.  
From a computation point of view, as $\nabla F_\mu$ can be efficiently evaluated through proximal mappings of both (weakly) convex components, it is possible to further incorporate quasi-Netwon updates \citep{byrd1994representations} or Nesterov-type accelerations \citep{nesterov1983method} for solving problem \eqref{eq: dc}. In addition, the proximal mappings of $h$ (or $\phi$ in GD) and $g$ can be computed simultaneously. Moreover, as we will see in the next section, the DME smoothing can be combined with the classical augmented Lagrangian to solve constrained programs \eqref{eq: lcdc}.}
\section{Linearly Constrained DC (LCDC) Optimization}\label{sec: LCDC}
In this section we consider the LCDC program \eqref{eq: lcdc}:
\begin{equation*}
	\min_{x\in \R^n} F(x) = \phi(x) - g(x)\quad \mathrm{s.t.} \quad Ax=b.
\end{equation*}
It is known that under some mild regularity condition at a local minimizer $x$ \cite[Theorem 8.15]{rockafellar2009variational}, there exists $\lambda \in \R^m$ such that
\begin{align}\label{eq: smooth lcdc stationary condition}
	0 \in \partial F(x) + A^\top \lambda, \ \text{and} \ Ax-b = 0.
\end{align}
Based on different structures of $\phi$, we propose two proximal ALM algorithms that will find an approximate solution to \eqref{eq: smooth lcdc stationary condition}. In particular, we present LCDC-ALM in Section \ref{sec: smooth LCDC} when $\phi$ is a smooth (possibly nonconvex) function, and composite LCDC-ALM in Section \ref{sec: nonsmooth LCDC} when $\phi$ is the sum of a smooth (possibly nonconvex) function and a nonsmooth convex function. 

\subsection{LCDC-ALM}\label{sec: smooth LCDC}
In this subsection, we make the following additional assumptions on the problem data. 
\begin{assumption}\label{assumption: smooth lcdc}~
\begin{enumerate}
	\item $\phi = f:\R^n\rightarrow \R$ is Lipschitz differentiable with modulus $L_f>0$, i.e., $h$ is the zero function.
    \item The matrix $A\in \R^{m\times n}$ is nonzero and $b \in \R^m$ belongs to the column space of $A$. 
	\item There exist {finite} $\overline{\mu}, \overline{\rho}>0$ such that for all $0 < \mu < \overline{\mu}$ and $\rho > \overline{\rho}$,
	\begin{align}\label{assumption: lower bounded}
		v(\mu,\rho) := \inf_{x,y\in \R^n} \left\{f(x) -g(y) + \frac{1}{2\mu}\|x-y\|^2 + \frac{\rho}{2}\|Ax-b\|^2\right\}> -\infty.
	\end{align}
	\end{enumerate}
\end{assumption}
We justify condition \eqref{assumption: lower bounded} in Lemma \ref{lemma: v lowerbounded}, which allows $g$ to be any Lipschitz convex function, or problem \eqref{eq: lcdc} to be any quadratic program that is regular in some sense. 
\begin{lemma}\label{lemma: v lowerbounded}
    Condition \eqref{assumption: lower bounded} can be satisfied under each of the following cases.
    \begin{enumerate}
        \item The function $f-g$ achieves a finite minimal value, and $g$ is $L_g$-Lipschitz continuous.
        \item The function $f$ is quadratic (possibly nonconvex), $g$ is $L_g$-Lipschitz continuous, and $x^\top \nabla^2 f x  > 0$ for all $x\in \R^n$ such that $Ax=0$. 
        \item The function $f$ is quadratic (possibly nonconvex), $g$ is a convex quadratic function, and  
        \begin{align}
            \min_{x\in \R^n}\{x^\top \nabla^2 f x~|~ Ax=0,\|x\|=1\}> \lambda_{\max}(\nabla^2 g). \notag 
        \end{align} 
    \end{enumerate}
    (Note that 2 and 3 coincide when $g$ is affine.) 
\end{lemma}
\proof{Proof.}
    See \ref{proof: v lowerbounded}. \Halmos
\endproof
Next we define approximate stationarity for problem \eqref{eq: lcdc} as follows.
\begin{definition}
	Given $\epsilon >0$, we say $x\in \R^n$ is an $\epsilon$-stationary point of \eqref{eq: lcdc} under Assumption \ref{assumption: smooth lcdc} if there exists $(\xi, y, \lambda)\in \R^n \times \R^n \times \R^m$ such that 
	\begin{align}\label{eq: smooth lcdc stationary condition approx}
		 \xi \in  \nabla f(x) - \partial g(y) + A^\top \lambda, \ \text{and} \ \max\{ \|\xi\|, \|x-y\|, \|Ax-b\|\} \leq \epsilon.
	\end{align}
\end{definition}
Before presenting the algorithm, we first introduce some notation. The classic augmented Lagrangian function is defined as
\begin{equation}
	L_{\rho}(x, \lambda) := f(x)-g(x) + \langle \lambda, Ax+b\rangle + \frac{\rho}{2}\|Ax-b\|^2,
\end{equation}
for some $\rho>0$. We also introduce a smoothed proximal augmented Lagrangian function,
\begin{equation}
	{\psi(x, z, \lambda) := \left(\tf_\rho(x,\lambda)+\frac{1}{2\mu}\|x-z\|^2\right) - M_{\mu g}(z)},
\end{equation}
{for some $\mu > 0$}, where $\tf_{\rho}(x,\lambda) := f(x) + \langle \lambda, Ax+b\rangle + \frac{\rho}{2}\|Ax-b\|^2$.
Furthermore, {the Moreau envelope $M_{\mu\tf_\rho}$ and proximal mapping $x_{\mu\tf_\rho}$ of $\tf_\rho$ are given below},
\begin{align*}
{M_{\mu \tilde{f}_\rho}(z, \lambda)=} \ {\min_{x\in \R^n}\left\{ \tilde{f}_\rho(x,\lambda)+\frac{1}{2\mu}\|x-z\|^2\right\}, \text{ and}} \ 
{x_{\mu\tilde{f}_\rho}(z, \lambda)=} \ {\argmin_{x\in \R^n}\left\{ \tilde{f}_\rho(x,\lambda)+\frac{1}{2\mu}\|x-z\|^2\right\}}.
\end{align*}
\begin{algorithm}[tbh]
	\caption{: LCDC-ALM} \label{alg: smooth lcdc-alm}
	\begin{algorithmic}[1]
		\STATE \textbf{Let} $0 < \mu < \min\{\overline{\mu},(L_f)^{-1}\}$, $0 < \beta < 2$, $\rho > \overline{\rho}$;
		\STATE \textbf{Initialize} {$x^{0}\in \R^n$,  $z^{0}\in \R^n$, $\lambda^{0}\in \R^m$};
		\FOR{$k = 0,1,\cdots$}
		\STATE \begin{align}\label{alg: smooth lcdc-alm x ineact}
 					\hspace{-3mm} x^{k+1} = \argmin_{x\in \R^n}\Big\{ & \langle \nabla f(x^k), x-x^k \rangle + \langle \lambda^k, Ax-b\rangle + \frac{\rho}{2}\|Ax-b\|^2 + \frac{1}{2\mu}\|x-z^k\|^2\Big \};
				\end{align}
		\STATE $x_{\mu g}(z^k) = \argmin_{x\in \R^n}\left\{g(x) + \frac{1}{2\mu}\|x-z^k\|^2\right\}$;
		\STATE $z^{k+1} = z^k + \beta (x^{k+1} - x_{\mu g}(z^k))$; \label{alg: smooth lcdc-alm z grad}
		\STATE $\lambda^{k+1} = \lambda^k + \rho (Ax^{k+1}-b)$; \label{alg: smooth lcdc-alm lmd grad}
		\ENDFOR
	\end{algorithmic}
\end{algorithm}
 
The LCDC-ALM is presented in Algorithm \ref{alg: smooth lcdc-alm}. To provide some intuition, notice that in analog to the unconstrained case, given some dual variable $\lambda\in \R^m$, the function 
\begin{equation*}
{\mathcal{L}_{\mu\rho}(z, \lambda)} := {M_{\mu \tf_\rho}(z,\lambda)} - M_{\mu g}(z)
\end{equation*}
is a smooth approximation of the augmented Lagrangian $L_{\rho}(x, \lambda)$ with gradient given by
\begin{align*}
	\nabla_z {\mathcal{L}_{\mu\rho}(z, \lambda)} = &\ \frac{1}{\mu}(z - x_{\mu\tf_\rho}(z, \lambda)) -  \frac{1}{\mu}(z -  x_{\mu g}(z)) = \frac{1}{\mu}(x_{\mu g}(z) - x_{\mu\tf_\rho}(z, \lambda)), {\text{ and}}\\
	{\nabla_\lambda \mathcal{L}_{\mu\rho}(z,\lambda)} = &\ {Ax_{\mu\tf_\rho}(z,\lambda)-b.}
\end{align*}
{Since the Moreau envelope $M_{\mu \tf_\rho}$ may still be difficult to evaluate, \eqref{alg: smooth lcdc-alm x ineact} replaces $f$ with its linearization and thus computes an inexact gradient of $M_{\mu\tf_\rho}$, which is equivalent to solving the following positive definite linear system}
\begin{align*}
		\left(\rho A^\top A + \frac{1}{\mu} I \right) x = \frac{1}{\mu}z^k + \rho A^\top b - A^\top \lambda^k-\nabla f(x^k).
\end{align*}
Consequently, line \ref{alg: smooth lcdc-alm z grad} and line \ref{alg: smooth lcdc-alm lmd grad} can be regarded as an {inexact} gradient descent of {$\mathcal{L}_{\mu\rho}(z,\lambda)$} in $z$ and gradient ascent of {$\mathcal{L}_{\mu\rho}(z,\lambda)$} in $\lambda$, respectively. As we will show later, this {inexact gradient descent on $z$} is sufficient to ensure convergence {of Algorithm \ref{alg: smooth lcdc-alm}}. Meanwhile, when {the Moreau envelope $M_{\mu\tf_\rho}$ is easy to compute,}
we can always replace \eqref{alg: smooth lcdc-alm x ineact} by its exact version, i.e., {$x^{k+1}=x_{\mu\tilde{f}_\rho}(z^k, \lambda^k)$,}
whose analysis is similar and omitted in this paper. 
We firstly establish some descent properties of $\psi$ in the next lemma.
\begin{lemma} \label{lemma: smooth lcdc-alm descent}
Suppose Assumptions \ref{assumption: dc} and \ref{assumption: smooth lcdc} hold. For all $k\in \N$, we have 
\begin{align*}
	\psi(x^k, z^k, \lambda^k)-\psi(x^{k+1}, z^{k+1}, \lambda^{k+1}) 
	\geq & \left(\frac{\mu^{-1}-L_f}{2} \right)\|x^{k+1}-x^k\|^2 + \frac{1}{\mu}\left( \frac{1}{\beta} - \frac{1}{2}\right) \|z^{k+1}-z^k\|^2 -\frac{1}{\rho}\|\lambda^{k+1}-\lambda^k\|^2.
\end{align*}
\end{lemma}
\proof{Proof.}
	See \ref{proof: smooth lcdc-alm descent}. \Halmos
\endproof

Next we show that, due to the smoothness of $f$, the difference of dual variables can be effectively bounded by the difference of primal variables. For notation consistency, denote $x^{-1}=x^0$ and $z^{-1} = x^0 +\mu(\nabla f(x^0)+A^\top \lambda^0).$
\begin{lemma}\label{lemma: bound dual by primal}
	Suppose Assumptions \ref{assumption: dc} and \ref{assumption: smooth lcdc} hold. For all $k\in \N$,
	\begin{align*}
	 \|\lambda^{k+1}-\lambda^k\| \leq & \frac{1}{\sigma_{\min}^+(A)} \left(\mu^{-1}\|x^{k+1}-x^k\| + L_f\|x^k-x^{k-1}\| + \mu^{-1}\|z^k-z^{k-1}\|\right).
	\end{align*}
\end{lemma}
\proof{Proof.}
	See \ref{proof: bound dual by primal}. \Halmos
\endproof

Lemmas \ref{lemma: smooth lcdc-alm descent} and \ref{lemma: bound dual by primal} suggest that the sequence defined by 
\begin{align}
	\Psi_k := \psi(x^k, z^k, \lambda^k) + \frac{\nu}{2}	\|x^k-x^{k-1}\|^2+\frac{\nu}{2}	\|z^k-z^{k-1}\|^2
\end{align}
is a decreasing sequence for some properly chosen $\nu>0$. To this end, recall $c_1$ and $c_2$ defined in \eqref{eq: c_12}, and we further let 
\begin{align}\label{eq: c_34}
	c_3 = & \frac{3\mu^{-2}}{\sigma^+_{\min}(AA^\top)}, \quad c_4 = \frac{3L_f^2}{\sigma^+_{\min}(AA^\top)}.
\end{align}
We assume that $\nu>0$ and $\rho>0$ are chosen such that  
\begin{align}\label{eq: kappa}
	\kappa_1 = c_1-\frac{{c}_3}{\rho}-\frac{{\nu}}{2} >0, \ \kappa_2 = c_2- \frac{{\nu}}{2}>0, \ \kappa_3 = \frac{{\nu}}{2}-\frac{{c}_4}{\rho}>0 ,\ \kappa_4= \frac{{\nu}}{2}-\frac{{c}_3}{\rho} >0.
\end{align}
This is always possible, for example, by first letting $\nu <2 \min \{c_1, c_2\}$ and then choosing $\rho > \max \left\{ \frac{c_3}{c_1-\nu/2},  \frac{2c_3}{\nu},  \frac{2c_4}{\nu}\right\}.$
\begin{lemma}\label{lemma: smooth lcdc-alm potential function}
	Suppose Assumptions \ref{assumption: dc} and \ref{assumption: smooth lcdc} hold. The following claims hold for LCDC-ALM.
	\begin{enumerate}
		\item For all $k\in\N$, we have 
			\begin{align} \label{eq: potentia descent}
			  	\Psi_k-{\Psi}_{k+1} \geq & \kappa_1 \|x^{k+1}-x^k\|^2 + \kappa_2 \|z^{k+1}-z^k\|^2 + \kappa_3 \|x^k-x^{k-1}\|^2 +\kappa_4 \|z^k-z^{k-1}\|^2. 
			  \end{align}
		\item $\Psi_k$ is bounded from below by $v(\mu, \rho)$ defined in \eqref{assumption: lower bounded} for all $k\in \N$.
		\item For any positive integer $K$, there exists an index $0\leq \bar{k}\leq K-1$ such that 
		\begin{align*}
			 \max\left \{\|x^{\bar{k}+1}-x^{\bar{k}}\|^2, \|z^{\bar{k}+1}-z^{\bar{k}}\|^2,\|x^{\bar{k}}-x^{\bar{k}-1}\|^2,\|z^{\bar{k}}-z^{\bar{k}-1}\|^2 \right\} \leq  \frac{\Psi_0 - v(\mu,\rho)}{\kappa_{\min} K },
		\end{align*}
		where $\kappa_{\min} =\min\{ \kappa_1, \kappa_2,  \kappa_3, \kappa_4\} >0$.
	\end{enumerate}
\end{lemma}
\proof{Proof.}
	See  \ref{proof: smooth lcdc-alm potential function} . \Halmos
\endproof
Utilizing Lemma \ref{lemma: smooth lcdc-alm descent}-Lemma \ref{lemma: smooth lcdc-alm potential function}, we are ready to present the convergence of LCDC-ALM.
\begin{theorem}\label{thm: smooth lcdc-alm}
	Suppose Assumptions \ref{assumption: dc} and \ref{assumption: smooth lcdc} hold. Let $\{(x^{k+1},z^{k+1},\lambda^{k+1})\}_{k\in \N}$ be the sequence generated by LCDC-ALM, and 
	denote $$\xi^{k+1} := \nabla f(x^{k+1})-\nabla f(x^k) + \mu^{-1}(x_{\mu g}(z^k) - x^{k+1}).$$ For any positive index $K$, there exists $0\leq \bar{k}\leq K-1$ such that $x^{\bar{k}+1}$, together with $$(\xi^{\bar{k}+1}, x_{\mu g}(z^{\bar{k}}), \lambda^{\bar{k}+1})\in \R^n\times \R^n \times \R^m,$$
	is an approximate stationary solution of \eqref{eq: lcdc} satisfying
	\begin{subequations}
	\begin{align}
		& \xi^{\bar{k}+1} \in \nabla f(x^{\bar{k}+1}) - \partial g(x_{\mu g}(z^{\bar{k}})) + A^\top \lambda^{k+1}, {\text{ and}} \label{eq: smooth lcdc-alm iter approx stationary a}\\
		& \max  \left \{ \|\xi^{\bar{k}+1}\|,  \|x_{\mu g}(z^{\bar{k}}) - x^{\bar{k}+1}\|, \|Ax^{\bar{k}+1}-b\|\right \} \leq C\left(\frac{{\Psi}_0 - {v(\mu,\rho)}}{\kappa_{\min} K}\right)^{1/2},  \label{eq: smooth lcdc-alm iter approx stationary b}
	\end{align}
	\end{subequations}
	where $v(\mu, \rho)$ is defined in \eqref{assumption: lower bounded}, $\kappa_{\min} =\min\{ \kappa_1, \kappa_2,  \kappa_3, \kappa_4\} >0$,
	\begin{align*}
		{\Psi}_0 = & \psi(x^0,z^0,\lambda^0)+ \frac{\nu}{2}\left\|x^0-z^0+\mu(\nabla f(x^0)+A^\top \lambda^0)\right\|^2, \ \text{and} \\ 
		C = & L_f + \frac{\mu^{-1}+1}{\beta} + \frac{2\sqrt{{c}_3/3}}{\rho}+ \frac{\sqrt{{c}_4/3}}{\rho}.
	\end{align*}
	That is, LCDC-ALM finds an $\epsilon$-approximate solution of \eqref{eq: lcdc} in the sense of \eqref{eq: smooth lcdc stationary condition approx} in no more than 
	\begin{align}\label{eq:  inexact iter complexity}
		\left \lceil  \frac{C^2 ({\Psi}_0 - v(\mu, \rho))}{\kappa_{\min}\epsilon^2} \right \rceil =\mathcal{O}(\epsilon^{-2})
	\end{align}
	iterations.
\end{theorem}
\proof{Proof.}
The optimality conditions of updates of $x^{k+1}$ and $x_{\mu g}(z^k)$ are given as
\begin{align*}
	0 =  \nabla f(x^{k}) + A^\top \lambda^{k+1} + \mu^{-1}(x^{k+1}-z^k),\text{ and} \;\, 0 \in  \partial g(x_{\mu g}(z^k)) + \mu^{-1}(x_{\mu g}(z^k)-z^k).
\end{align*}
As a result, we have
\begin{align}
	\xi^{k+1} = &  \nabla f(x^{k+1})-\nabla f(x^k) + \mu^{-1}(x_{\mu g}(z^k) - x^{k+1}) 
	 \in    \nabla f(x^{k+1}) -\partial g(x_{\mu g}(z^k))+ A^\top \lambda^{k+1}.\label{eq:smoothLCDC-ALM xi}
\end{align}
Since $x^{k+1}-x_{\mu g}(z^k) = \frac{1}{\beta}(z^{k+1}-z^k)$ and $Ax^{k+1}-b = \frac{1}{\rho}(\lambda^{k+1}-\lambda^k)$, we have
\begin{align*}
	 &\max  \left \{\|\xi^{k+1}\|, \; \|x_{\mu g}(z^k) - x^{k+1}\|, \; \|Ax^{k+1}-b\|\right \}\\
	\leq &  \max  \Big \{ L_f \|x^{k+1}-x^k\| + \frac{1}{\mu\beta}\|z^{k+1}-z^k\|, \;\; \frac{1}{\beta}\|z^{k+1}-z^k\|, \\
		 &\quad \quad \quad  \frac{\sqrt{{c}_3/3}}{\rho} \|x^{k+1}-x^k\|  + \frac{\sqrt{{c}_4/3}}{\rho} \|x^{k}-x^{k-1}\| + \frac{\sqrt{{c}_3/3}}{\rho} \|z^{k}-z^{k-1}\| \Big\}\\
	 \leq & \ C \max \left\{ \|x^{k+1}-x^k\|,  \|x^{k}-x^{k-1}\|, \|z^{k+1} - z^{k}\|, \|z^{k} - z^{k-1}\| \right\},
\end{align*}
where the {first} inequality is due to Lemma \ref{lemma: bound dual by primal}. Invoking Lemma \ref{lemma: smooth lcdc-alm potential function}, there exists $0\leq \bar{k}\leq K-1$ such that 
\begin{align*}
	\max  \left \{ \|\xi^{\bar{k}+1}\|,  \|x_{\mu g}(z^{\bar{k}}) - x^{\bar{k}+1}\|, \|Ax^{\bar{k}+1}-b\|\right \} \leq  C\left(\frac{{\Psi}_0 - {v(\mu,\rho)}}{\kappa_{\min} K}\right)^{1/2}. 
\end{align*}
This proves \eqref{eq: smooth lcdc-alm iter approx stationary b} and \eqref{eq:  inexact iter complexity} and substituting $\bar{k}$ into \eqref{eq:smoothLCDC-ALM xi} proves \eqref{eq: smooth lcdc-alm iter approx stationary a}. \Halmos
\endproof
\subsection{Composite LCDC-ALM}\label{sec: nonsmooth LCDC}
In this subsection, we consider a more challenging setup of problem \eqref{eq: lcdc}, where we allow $\phi$ to have a nonsmooth component $h$. In particular, we make some additional assumptions on $h$.
\begin{assumption}\label{assumption: nonsmooth lcdc}
	In addition to being proper, closed, and convex, the function $h:\R^n\rightarrow \overline{\R}$ is $L_h$-Lipschitz {continuous} over its effective domain $\mathcal{H}=\{x\in \R^n~|~h(x)<+\infty\}$, {which is nonempty, convex, and compact, i.e.,} there exists $0 < D_\mathcal{H} < +\infty$ such that $	\|x^1-x^2\| \leq  D_\mathcal{H}$ for all $x^1, x^2 \in \mathcal{H}$. Moreover, there exists $\bar{x} \in \mathrm{int} ~\mathcal{H}$ (interior of $\mathcal{H}$) such that $A\bar{x}=b$;
\end{assumption} 
By Assumption \ref{assumption: nonsmooth lcdc}, we can further define positive constants 
{\begin{align}
	M_{\nabla f} := \max_{x\in \mathcal{H}} \|\nabla f(x)\|, \ \bar{d} := \mathrm{dist} (\bar{x}, \partial \mathcal{H}), \ \text{and} \ M_{\partial g} := \sup_{\xi_g \in \partial g(\mathcal{H})} \|\xi_g\| .
\end{align}
The constant $M_{\nabla f}$ is well-defined due to the continuity of $\nabla f$ over compact $\mathcal{H}$; $\partial \mathcal{H}$ is the boundary of $\mathcal{H}$, and hence $\bar{d}$ is strictly positive. The constant $M_{\partial g}$ is finite due to \cite[Theorem 24.7]{rockafellar1970convex}.} Next we define an approximate stationary solution under Assumptions \ref{assumption: dc} and \ref{assumption: nonsmooth lcdc}.
\begin{definition}
	Given $\epsilon >0$, we say that $x\in \R^n$ is an $\epsilon$-stationary point of \eqref{eq: lcdc} under Assumptions \ref{assumption: dc} and \ref{assumption: nonsmooth lcdc} if there exists $(\xi, y, \lambda)\in \R^n \times \R^n \times \R^m$ such that 
	\begin{align}\label{eq: Nonsmooth lcdc stationary condition approx}
		\xi \in  \nabla f(x) +\partial h(x) - \partial g(y) + A^\top \lambda, \ \text{and} \ \max\{ \|\xi\|, \|x-y\|, \|Ax-b\|\} \leq \epsilon.
	\end{align}
\end{definition}
We present the composite LCDC-ALM in Algorithm \ref{alg: ns-lcdc-alm}.
\begin{algorithm}[tbh!]
	\caption{: Composite LCDC-ALM} \label{alg: ns-lcdc-alm}
	\begin{algorithmic}[1]
		\STATE \textbf{Input} $0 < \mu <  L_f^{-1}$, $0 < \beta \leq 1$, $\rho > 0$, and $\{\epsilon_{k+1}\}_{k\in \N} \subset [0,1]$ with $E := \sum_{k=1}^{\infty}\epsilon^2_{k} < +\infty$;
		\STATE \textbf{Initialize} {$x^{0}\in \mathcal{H}$,  $z^{0}\in  \mathcal{H}$, $\lambda^{0}\in \mathrm{Im}(A)$;}
		\FOR{$k = 0,1,\cdots$}
        \STATE evaluate $\xi_g^k \in \partial g(x^k)$ and find $(x^{k+1}, \zeta^{k+1}) \in \mathcal{H} \times \R^n$ with $\|\zeta^{k+1}\|\leq \epsilon_{k+1}$ such that 
		\begin{align}\label{eq: composite-lcdc-primal-update}
			\zeta^{k+1} \in \nabla f(x^k) - \xi_g^k + \partial h(x^{k+1}) + A^\top \lambda^k + \rho A^\top (Ax^{k+1}-b) + \frac{1}{\mu}(x^{k+1}-z^k);
		\end{align}s
        \vspace{-0.3in}
		\STATE $z^{k+1} = z^k + \beta (x^{k+1} - z^k)$;
		\STATE $\lambda^{k+1} = \lambda^k + \rho (Ax^{k+1}-b)$;
		\ENDFOR
	\end{algorithmic}
\end{algorithm}
\begin{remark}
Due to the additional nonsmooth function $h$, we cannot effectively control the dual difference using primal iterates as in Lemma \ref{lemma: bound dual by primal}. We simply replace $-g$ by its linearization at $x^k$, and then smooth the linearized AL function: condition \eqref{eq: composite-lcdc-primal-update} says that $x^{k+1}$ is an approximate solution of the following problem: 
\begin{align}\label{eq: composite update in opt form}
    \min_{x\in \R^n} \langle \nabla f(x^k)-\xi_g^k , x-x^k \rangle + h(x) + \langle \lambda^k, Ax-b\rangle + \frac{\rho}{2}\|Ax-b\|^2 + \frac{1}{2\mu}\|x-z^k\|^2.
\end{align}
Though not directly using the {DME} smooth approximation \eqref{eq: smooth approx}, Algorithm \ref{alg: ns-lcdc-alm} still relies on the smoothing property of the Moreau envelope. It can be regarded as a primal-dual version of the linearized pDCA, and is closely related to the proximal ALM proposed in \cite{zhang2020global,zhang2020proximal} and LiMEAL in \cite{zeng2021moreau}.
\end{remark}
Throughout the analysis, we assume that $\{\epsilon_{k+1}\}_{k\in \N} \subset [0,1]$ with $E := \sum_{k=1}^{\infty}\epsilon^2_{k} < +\infty$. We utilize the following proximal augmented Lagrangian function:
\begin{equation}
	P(x, z, \lambda) := f(x)+h(x)-g(x)+\langle \lambda, Ax-b\rangle + \frac{\rho}{2}\|Ax-b\|^2 + \frac{1}{2\mu}\|x-z\|^2,
\end{equation}
which will serve as a potential function as shown in the following lemma. 
\begin{lemma} \label{lemma: ns_lcdc_descent}
Suppose {Assumptions \ref{assumption: dc} and \ref{assumption: nonsmooth lcdc} hold}. For all $k\in \N$, we have 
\begin{align}
	& P(x^k, z^k, \lambda^k)-P(x^{k+1}, z^{k+1}, \lambda^{k+1}) \notag \\
	\geq & {\left( \frac{\mu^{-1}-2L_f}{4}\right)} \|x^{k+1}-x^k\|^2 + \frac{1}{2\beta\mu} \|z^{k+1}-z^k\|^2 -\frac{1}{\rho}\|\lambda^{k+1}-\lambda^k\|^2 -\mu \epsilon_{k+1}^2.  \label{eq: ns_descent}
\end{align}
\end{lemma}
\proof{Proof.}
	See \ref{proof: ns_lcdc_descent}. \Halmos
\endproof

Next we show that $P(x^1,z^1,\lambda^1)$ is bounded from above by a constant independent of $\rho$.
\begin{lemma}\label{lemma: bound init val}
	Suppose Assumptions \ref{assumption: dc} and \ref{assumption: nonsmooth lcdc} hold. Then $P(x^1,z^1,\lambda^1)\leq \overline{P}$, where 
	\begin{align}
	  \overline{P} := &  \max_{x\in \mathcal{H}} \{f(x)+h(x)-g(x)\} + 3(L_h+ M_{\nabla f} + M_{\partial g})D_{\mathcal{H}} + \notag \\
    & { \frac{L_f+ 6\mu^{-1}}{2} D_{\mathcal{H}}^2  + \frac{3L_f^{-1}}{2} } +2\|\lambda^{0}\|\max_{x\in \mathcal{H}}\|Ax-b\|. \label{eq: upper bound independent of rho}
	\end{align}
\end{lemma}
\proof{Proof.}
	See  \ref{proof: bound init val}. \Halmos
\endproof

Now we show that the sequence $\{\lambda^k\}_{k\in \N}$ stays bounded. 
\begin{lemma}\label{lemma: ns dual bd}
Suppose {Assumptions \ref{assumption: dc} and \ref{assumption: nonsmooth lcdc} hold}. For all $k\in \N$, we have 
\begin{align}\label{eq: ns dual bound}
	\|\lambda^k\| \leq \Lambda :=  \max\left\{\|{\lambda^{0}}\|,\, \frac{2D_\mathcal{H}}{\bar{d} {\sigma_{\min}^+(A)}}\left(M_{\nabla f}+ M_{\partial g}+\frac{D_{\mathcal{H}}}{\mu}+L_h {+1}\right)\right\}.
\end{align}
\end{lemma}
\proof{Proof.}
	See  \ref{proof: ns dual bd}. \Halmos
\endproof

In the next lemma, we show that $P(x^k, z^k, \lambda^k)$ is bounded from below due to the boundedness of $\lambda^k$ established in Lemma \ref{lemma: ns dual bd}, which further allows us to bound the difference of consecutive primal iterates. 
\begin{lemma}\label{lemma: ns_lcdc}
Suppose {Assumptions \ref{assumption: dc} and \ref{assumption: nonsmooth lcdc} hold}. The following statements hold.
\begin{enumerate}
	\item Recall $\Lambda$ from Lemma \ref{lemma: ns dual bd}. For all $k\in \N$, we have
 		\begin{align}\label{eq: ns lb p}
			P(x^{k}, z^k, \lambda^k) \geq \underline{P}:= \min_{x\in \R^n} \left\{f(x)+h(x)-g(x)\right\}- \Lambda \max_{x\in \mathcal{H}}\|Ax-b\|	> -\infty. 
		\end{align}
	\item Recall $\overline{P}$ in \eqref{eq: upper bound independent of rho}, {the constant $E = \sum_{k=1}^\infty \epsilon_k^2 < +\infty$,} and define 
    \begin{align}
         {\eta := \min\{\frac{1}{4}(\mu^{-1}-2L_f), (2\mu \beta)^{-1}\}. }
    \end{align}
	For any positive $K\in \N$, there exists ${1}\leq \bar{k}\leq {K} $ such that 
	\begin{align}
		\max \left\{ \|x^{\bar{k}+1} -x^{\bar{k}}\|^2, \|z^{\bar{k}+1} -z^{\bar{k}}\|^2\right\} \leq  \frac{{\overline{P}}-\underline{P} {+\mu E}}{\eta K} + \frac{8 \Lambda^2}{\eta \rho}.
	\end{align} 
\end{enumerate}	
\end{lemma}
\proof{Proof.}
	See  \ref{proof: ns_lcdc}. \Halmos
\endproof

Now we are ready to present the convergence of composite LCDC-ALM.
\begin{theorem}\label{thm: ns LCDC-ALM}
	Suppose {Assumptions \ref{assumption: dc} and \ref{assumption: nonsmooth lcdc} hold}. Let $\{(x^{k+1},z^{k+1},\lambda^{k+1})\}_{k\in \N}$ be the sequence generated by composite LCDC-ALM, and define $\xi^{k+1} := {\zeta^{k+1}}+ \nabla f(x^{k+1})-\nabla f(x^k) + \mu^{-1}(z^k - x^{k+1}).$ For any positive $K\in \N$, there exists {$1\leq \bar{k}\leq K$} such that $x^{\bar{k}+1}$, together with $(\xi^{\bar{k}+1}, x^{\bar{k}}, \lambda^{\bar{k}+1})\in \R^n\times \R^n \times \R^m,$
	is an approximate stationary solution of problem \eqref{eq: lcdc} satisfying
	\begin{subequations}\label{eq: ns lcdc-rate}
	\begin{align}
		& \xi^{\bar{k}+1} \in \nabla f(x^{\bar{k}+1}) + \partial h(x^{\bar{k}+1}) - \partial g(x^{\bar{k}}) + A^\top \lambda^{\bar{k}+1}, \label{eq: composite kkt1}\\
		& \max  \left \{ \|\xi^{\bar{k}+1}\|,  \|x^{\bar{k}+1}-x^{\bar{k}}\|\right \} \leq  \gamma \left( \frac{\overline{P}-\underline{P} {+\mu E}}{\eta K} + \frac{8 \Lambda^2}{\eta\rho}\right)^{1/2} {+ \epsilon_{\bar{k}+1}},\label{eq: ns lscd-rate dual res}\\
		& \|Ax^{\bar{k}+1}-b\| \leq  \frac{2\Lambda}{\rho}, 
	\end{align}
	\end{subequations}
	where {$\overline{P}$ is given in \eqref{eq: upper bound independent of rho}}, $\underline{P}$ is defined in \eqref{eq: ns lb p}, $\Lambda$ is defined in \eqref{eq: ns dual bound}, {$E=\sum_{k=1}^{\infty} \epsilon^2_{k}$}, $\gamma = L_f + 1/(\mu\beta)+1$, and {$\eta = \min\{(\mu^{-1}-2L_f)/4, 1/(2\mu \beta)\}$.}
	In other words, {further suppose $\epsilon >0$ such that $0\leq \epsilon_k \leq \epsilon/2$ for all positive integer $k$.} If we choose 
	\begin{equation}\label{eq: rho bd}
		\rho \geq \max\left\{ \frac{2\Lambda}{\epsilon},  \frac{{64}\Lambda^2\gamma^2}{\eta \epsilon^2} \right\},  
	\end{equation}
	then composite LCDC-ALM finds an $\epsilon$-approximate solution of \eqref{eq: lcdc} in the sense of \eqref{eq: Nonsmooth lcdc stationary condition approx} in no more than $K+1$ iterations, where 
	\begin{equation}\label{eq: iter complex}
		K \leq \left \lceil  \frac{{8}\gamma^2({\overline{P}}-\underline{P} {+\mu E})}{\eta \epsilon^2} \right \rceil =\mathcal{O}(\epsilon^{-2})
	\end{equation}
    iterations.
\end{theorem}
\proof{Proof.}
	Let $\bar{k}$ be the index given in Lemma \ref{lemma: ns_lcdc}. The optimality condition of $x^{\bar{k}+1}$ gives \eqref{eq: composite kkt1} immediately. 
	 It holds that $\max  \left \{ \|\xi^{\bar{k}+1}\|,  \|x^{\bar{k}+1}-x^{\bar{k}}\|\right \} $ is bounded from above by
	\begin{align}
	   & \max \left\{{\|\zeta^{\bar{k}+1}\| + } \| \nabla f(x^{\bar{k}+1})-\nabla f(x^{\bar{k}})\|	 + \|\mu^{-1}(z^{\bar{k}} - x^{\bar{k}+1})\|,  \|x^{\bar{k}+1}-x^{\bar{k}}\| \right\} \notag \\
		\leq & \left(L_f+ \frac{1}{\mu\beta}+1\right)\max \left\{ \|x^{\bar{k}+1} -x^{\bar{k}}\|, \|z^{\bar{k}+1} -z^{\bar{k}}\|\right\}  {+ \epsilon_{\bar{k}+1} } \leq   \gamma \left( \frac{{\overline{P}}-\underline{P} {+\mu E}}{\eta K} + \frac{8 \Lambda^2}{\eta\rho}\right)^{1/2}  {+ \epsilon_{\bar{k}+1} } \notag;
	\end{align}
	in addition, 
	\begin{align*}
		\|Ax^{\bar{k}+1}-b\|= \frac{\|\lambda^{\bar{k}+1}-\lambda^{\bar{k}}\|}{\rho} \leq \frac{\|\lambda^{\bar{k}+1}\|+\|\lambda^{\bar{k}}\|}{\rho}\leq \frac{2\Lambda}{\rho}.
	\end{align*}
	This proves \eqref{eq: ns lcdc-rate}. Finally recall that $\epsilon_{\bar{k}+1} \leq \epsilon/2$. It is straightforward to verify that the claimed lower bound of $\rho$ in \eqref{eq: rho bd} and upper bound of $K$ in \eqref{eq: iter complex} together ensure $\|Ax^{\bar{k}+1}-b\|\leq \epsilon$ and
	$$ \max  \left \{ \|\xi^{\bar{k}+1}\|,  \|x^{\bar{k}+1}-x^{\bar{k}}\|\right \} \leq \gamma \left(\frac{\epsilon^2}{{8}\gamma^2} + \frac{\epsilon^2}{{8}\gamma^2}\right)^{1/2} {+\frac{\epsilon}{2}}= \epsilon.$$
	This completes the proof. \Halmos
\endproof

{Before ending this section, we briefly discuss the \textit{first-order} iteration complexity of composite LCDC-ALM. In view of \eqref{eq: composite update in opt form}, the major computational task in each iteration of composite LCDC-ALM is to solve a strongly convex optimization problem of the form $\min_{x\in \R^n} \psi(x) + h(x)$, where, for a fixed index $k$,  $\psi(x) = \langle f(x^k) - \xi_g^k, x -x^k\rangle + \langle \lambda^k, Ax -b\rangle + \frac{\rho}{2}\|Ax -b\|^2 + \frac{1}{2\mu}\|x -x^k\|^2$ for $x\in \R^n$. It is easy to verify that  $\psi$ is $m_\psi$-strongly convex and has a $L_{\psi}$-Lipschitz gradient, where $m_{\psi} = \mu^{-1}$ and $L_{\psi} = \rho\|A\|^2 + \mu^{-1}$. We take the accelerated proximal gradient (APG) method adopted by \cite{li2021rate} as an example of an optimal first-order method for solving this strongly convex program. It finds $(x^{k+1}, \zeta^{k+1}) \in \mathcal{H} \times \R^n$ that satisfies $\|\zeta^{k+1}\|\leq \epsilon_{k+1}$ and \eqref{eq: composite-lcdc-primal-update} in no more than
\begin{align*}
    T_{k+1} := \left \lceil \sqrt{\frac{L_{\psi}}{m_{\psi}}} \log \frac{64 L_\psi^3 D^2_{\mathcal{H}}}{\epsilon_{k+1}^2 \mu}   + 1 \right\rceil
\end{align*}
iterations \citep[Lemma 1]{li2021rate}. We aim to bound the summation $\sum_{k=0}^{K} T_{k+1}$, where $K = \mathcal{O}(\epsilon^{-2})$ by \eqref{eq: iter complex}. Now let us suppose $\rho \geq \mu^{-1}/\|A\|^2$ so that $L_\phi \leq 2\rho \|A\|^2$, and set $\epsilon_{k} = \epsilon/(2k)$ for all positive integer $k$ and some $\epsilon\in (0,2]$ (we consider $\mu$ as a constant independent of $\epsilon$).  We then have 
\begin{align*} 
    \sqrt{L_\psi/ m_\psi }  = \sqrt{L_{\psi} \mu }  \leq  \sqrt{ 2 \rho  \|A\|^2 L_f^{-1} }, \quad  \frac{64 L_\psi^3 D^2_{\mathcal{H}}}{\epsilon_{k+1}^2 \mu} = \frac{256 L_\psi^3 D^2_{\mathcal{H}} (k+1)^2}{\epsilon^2 \mu} \leq \frac{2048 D^2_{\mathcal{H}}\|A\|^8 \rho^4 (k+1)^2}{\epsilon^2}.
\end{align*}
By \eqref{eq: rho bd}, it suffices to choose $\rho = \mathcal{O}(\epsilon^{-2})$; by the above, for $k\leq K = \mathcal{O}(\epsilon^{-2})$, we can bound $T_{k+1}$ by
\begin{align*}
   T_{k+1} \leq \left \lceil  \sqrt{2\mathcal{O}(\epsilon^{-2}) \|A\|^2 L_f^{-1}} \log  \left( 2048D^2_{\mathcal{H}} \|A\|^8  \mathcal{O}(\epsilon^{-14})\right) +1 \right \rceil = \tilde{\mathcal{O}}(\epsilon^{-1}),
\end{align*}
where the $\tilde{\mathcal{O}}$ notation hides logarithmic dependency on $\epsilon$.  Consequently, the overall first-order iteration complexity of composite LCDC-ALM is given by 
\begin{align*}
    \sum_{k=0}^{K} T_{k+1}  \leq (K+1) \tilde{\mathcal{O}}(\epsilon^{-1}) = \tilde{\mathcal{O}}(\epsilon^{-3}).
\end{align*}
We note that this complexity upper bound matches existing results of ALM-based algorithms \citep{melo2020iteration, melo2020almfullstep} for affine-constrained weakly convex minimization. An interesting future direction is to explore the possibility of reducing the iteration complexity to $\tilde{\mathcal{O}}(\epsilon^{-2.5})$, which has been achieved when the concave component $-g$ is absent \citep{LiXu2020almv2}.}
\section{Numerical Experiments}\label{sec: numerical}
In this section, we present experiments to demonstrate the performance of proposed algorithms.
\subsection{Unconstrained DC Regularized Least Square Problem}
In this subsection, we consider the $\ell_{1-2}$ regularized least squares problem:
\begin{align}\label{eq: dcls}
	\min_{x\in \R^n} F(x) = \frac{1}{2}\|Cx-d\|^2 + \varrho \|x\|_1 -  \varrho \|x\|,
\end{align}
and compare the Inexact GD method (Algorithm \ref{alg: inexact gd on F_mu}) with the pDCAe proposed in \cite{wen2018proximal}. We generate data as described in \cite[Section 5]{wen2018proximal}: first create $C\in \R^{m\times n}$ with standard Gaussian entries, and then normalize its columns {to unit length}; generate $\hat{x}\in \R^n$ such that $\|\hat{x}\|_0 = s$ with nonzero Gaussian entries, and finally set $d = C \hat{x} + 0.01 \xi$, where $\xi\in \R^m$ has standard Gaussian entries. In our experiments, we choose $\beta = 1$ and {$\mu = 1/\|A\|^2$}. For pDCAe, the extrapolation parameters are chosen according to Section 3 of their paper, which are also popular choices used in FISTA \citep{beck2009fast}. Both algorithms are provided with the same initial point and will terminate if $\|x^{k+1}-y^k\|/\max\{1,\|x^{k+1}\|\}\leq 10^{-5}$, where $y^k=x_{\mu g}(z^k)$ in Inexact GD and $y^k=x^k$ in pDCAe.

\begin{table}[h!]
\caption{Comparison of Inexact GD method (Algorithm \ref{alg: inexact gd on F_mu}) and pDCAe \cite{wen2018proximal}}\label{table: DCLS}
\begin{center}
\begin{tabular}{lccccccc}
\hline
  &      & \multicolumn{2}{c}{Avg. Iteration} & \multicolumn{2}{c} {Avg. Time (s)}        & \multicolumn{2}{c}{Avg. Time{/Iter ($10^{-2}$s)}}\\
$i$ & $\varrho$  & Inexact GD             & pDCAe  & Inexact GD             & pDCAe          & Inexact GD        & pDCAe             \\
\hline
  & 1    & \textbf{124}    & 1920   & \textbf{2.00}  & 4.29           & 1.60      & \textbf{0.22}     \\
1 & 0.1  & \textbf{174}    & 2988   & \textbf{2.69}  & 6.34           & 1.54      & \textbf{0.21}     \\
  & 0.01 & \textbf{1079}   & 1541   & 16.81          & \textbf{3.27}  & 1.56      & \textbf{0.21}     \\
\hline
  & 1    & \textbf{107}    & 2273   & \textbf{7.90}  & 22.28          & 7.35      & \textbf{0.98}     \\
2 & 0.1  & \textbf{194}    & 3005   & \textbf{14.06} & 29.35          & 7.26      & \textbf{0.98}     \\
  & 0.01 & \textbf{1077}   & 1572   & 75.56          & \textbf{15.38} & 7.02      & \textbf{0.98}     \\
\hline
  & 1    & \textbf{104}    & 2429   & \textbf{15.44} & 60.19          & 14.79     & \textbf{2.48}     \\
3 & 0.1  & \textbf{196}    & 3005   & \textbf{28.64} & 72.32          & 14.61     & \textbf{2.41}     \\
  & 0.01 & \textbf{1052}   & 1467   & 159.98         & \textbf{37.11} & 15.21     & \textbf{2.53}    \\
\hline
\end{tabular}
\end{center}
\end{table}

For each $(m,n,s) = (720i, 2560i, 80i)$ where $i\in [3]$ and $\varrho \in \{1.0, 0.1, 0.01\}$, we generate five instances and report the average iteration number and wall clock time in Table \ref{table: DCLS}. Since both algorithms terminate with the same objective value in most cases, the average objective is the same up to four decimal places and therefore omitted from comparison. Inexact GD requires fewer iterations than pDCAe, and terminates faster in wall clock time when $\varrho \in \{1.0, 0.1\}$. We give two possible explanations for the slowness of Inexact GD when $\varrho = 0.01$: 1) since in general a proximal evaluation is more expensive than a subgradient evaluation, the per-iteration time of Inexact GD can be higher than pDCAe, {while a proper parallel implement might be the remedy,} and 2) when a smaller $\varrho$ is used, the concave component $-\|x\|$ in problem \eqref{eq: dcls} tends to vanish, and pDCAe behaves more like FISTA applied to convex composite optimization problems, whose efficiency is well-recognized.

\subsection{Constrained DC Regularized Least Square Problem}
In this subsection, we consider the DC regularized least squares problem with affine constraints:
\begin{align}\label{eq: c-dcls}
	\min_{x\in \R^n} F(x) = \frac{1}{2}\|Cx-d\|^2 +\delta_{\{x:\|x\|_1\leq M\}}(x)-  \varrho \|x\|\quad \mathrm{s.t.}\quad Ax=b,
\end{align} 
where we  explicitly require $\|x\|_1\leq M$ instead of penalizing $\|x\|_1$ as in \eqref{eq: dcls}. We generate $C\in \R^{m\times n}$ and $d\in \R^m$ the same way as in the previous subsection, $A\in \R^{m\times n}$ with standard Gaussian entries, and $b = A\tilde{x}$, where each component of $\tilde{x}$ is uniformly sampled from $[-M/(2n), M/(2n)]$. We compare the composite LCDC-ALM (Algorithm \ref{alg: ns-lcdc-alm}) with GD \eqref{eq: GD}, {DCA, and pDCA \eqref{eq: pDCA}} according to the following problem decomposition: in GD, {DCA, and pDCA}, let $\phi(x)=\frac{1}{2}\|Cx-d\|^2 +\delta_{\{x:Ax=b,\|x\|_1\leq M\}}(x)$; in composite LCDC-ALM, let $f(x) =\frac{1}{2}\|Cx-d\|^2$ and $h(x)=\delta_{\{x:\|x\|_1\leq M\}}(x)$; in both algorithms, we set $g(x)=\varrho \|x\|$. 

The infeasibility of problem \eqref{eq: c-dcls} is measured by $\|Ax^{k+1}-b\|$ in {composite} LCDC-ALM, while the constraints $Ax=b$ are always satisfied in the GD, {DCA, and pDCA} as part of $\phi$. Therefore we compare the {four} algorithms as follows: we first run {composite} LCDC-ALM until $\|Ax^{k}-b\|\leq 10^{-5}$ and $\lvert F(x^{k})-F(x^{k-1})\rvert/\lvert F(x^{k})\rvert \leq 10^{-3}$; then we execute {GD, DCA, and pDCA} for the same number of iterations. All nonlinear subproblems are solved by IPOPT with linear solver MA57.  In our experiments, we use $\varrho = 1.0$, $M=2.0$ for problem \eqref{eq: c-dcls}, and set $\beta = 0.1$, $\mu=1/L_f$ {(this is also the proximal coefficient used in pDCA)}. For each $(m,n,s) = (50i, 200i, 10i)$ where $i\in [5]$, we generate five instances and report the average iteration, objective value, and running time in Table \ref{table: compare gd-alm}. {The objective function is evaluated at the last iterate for all algorithms. We observed that all algorithms converge to solutions with similar quality, while DCA and pDCA seem to achieve slightly better objective values. Meanwhile, the composite LCDC-ALM has a clear advantage in solution time. }

\begin{table}[h!]
\caption{{Comparison of Composite LCDC-ALM (ALM), GD, DCA, and pDCA}}\label{table: compare gd-alm}
\begin{center}
{\begin{tabular}{cc|cccc|cccc}
\hline
  &            & \multicolumn{4}{c|}{Avg. Objective}         & \multicolumn{4}{c}{Avg. Time (s)}      \\
$i$ & Iteration & ALM     & GD      & DCA     & pDCA    & ALM    & GD     & DCA    & pDCA   \\
\hline
1 & 104        & 2.7564  & 2.7574  & 2.7517  & \textbf{2.7490}  & \textbf{3.30}   & 4.55   & 4.89   & 4.75   \\
2 & 113        & 8.8506  & 8.8514  & \textbf{8.8494}  & 8.8500  & \textbf{16.85}  & 22.25  & 24.28  & 22.61  \\
3 & 115        & 10.4611 & 10.4612 & 10.4608 & \textbf{10.4606} & \textbf{42.25}  & 70.28  & 80.10  & 71.95  \\
4 & 94         & 17.4131 & 17.4134 & 17.4131 & \textbf{17.4127} & \textbf{73.71}  & 114.13 & 127.00 & 111.19 \\
5 & 89         & 24.1486 & 24.1489 & \textbf{24.1483} & 24.1483 &\textbf{ 120.49} & 170.40 & 193.03 & 176.45 \\
\hline
\end{tabular}}
\end{center}
\end{table}

\subsection{Linearly Constrained Nonconvex Quadratic Program}
In this subsection, we consider the linearly constrained nonconvex quadratic program:
\begin{align}
	\min_{x\in \R^n} F(x) = f(x)-g(x) \quad \mathrm{s.t.} \quad Ax=b,
\end{align}
where $f(x) = \frac{1}{2} x^\top Qx + q^\top x$ and $g(x) = \frac{1}{2}x^\top G x$ for some $Q, G\in \mathbb{S}_+^{n \times n}$ and $q \in \R^n$. 
In particular, we would like to investigate how the choice of $\beta \in (0,2)$ affects the convergence of LCDC-ALM (Algorithm \ref{alg: smooth lcdc-alm}), and then compare its performance with the proximal ALM proposed in \cite{zhang2020global}. We first generate $A \in \R^{m \times n}$,  $q\in \R^n$, and $\hat{x} \in \R^n$ with standard Gaussian entries, and set $b = A\hat{x}$. Suppose for simplicity $\mathrm{rank}(A) = m$. Let $\{v_1,\cdots, v_{n-m}\}$ and $\{u_1,\cdots, u_m\}$ be an orthonormal basis of the null space of $A$ and the column space of $A^\top$, respectively; further denote $V_+ = \{v_1,\cdots, v_m,$ $u_1,\cdots, u_{\lfloor m/2 \rfloor}\}$ and $V_- = \{ u_{\lfloor m/2 \rfloor +1},\cdots, u_m\}$. Then we let $Q = \sum_{v\in V_+} a(v) vv^\top$ and $G = \sum_{u\in V_-} b(u) uu^\top$, where each $a(v)$ is uniformly sampled from $[0,10]$, and each $b(u)$ is uniformly sampled from $[0,50]$. The construction ensures $\inf_x \{F(x): Ax=b\} > -\infty$ and Assumption \ref{assumption: smooth lcdc} is satisfied. In our experiments, we use $m = 200$ and $n=500$; the spectrum of the generated matrix $Q-G$ lies between $[-49.74, 9.97]$.

\begin{figure*}[h!]
\caption{Infeasibility and Objective Trajectory of LCDC-ALM (Algorithm \ref{alg: smooth lcdc-alm})}
\center{
\label{fig: compare beta}  
	\begin{tabular}{@{}c@{}}
    	\includegraphics[width=.4\linewidth]{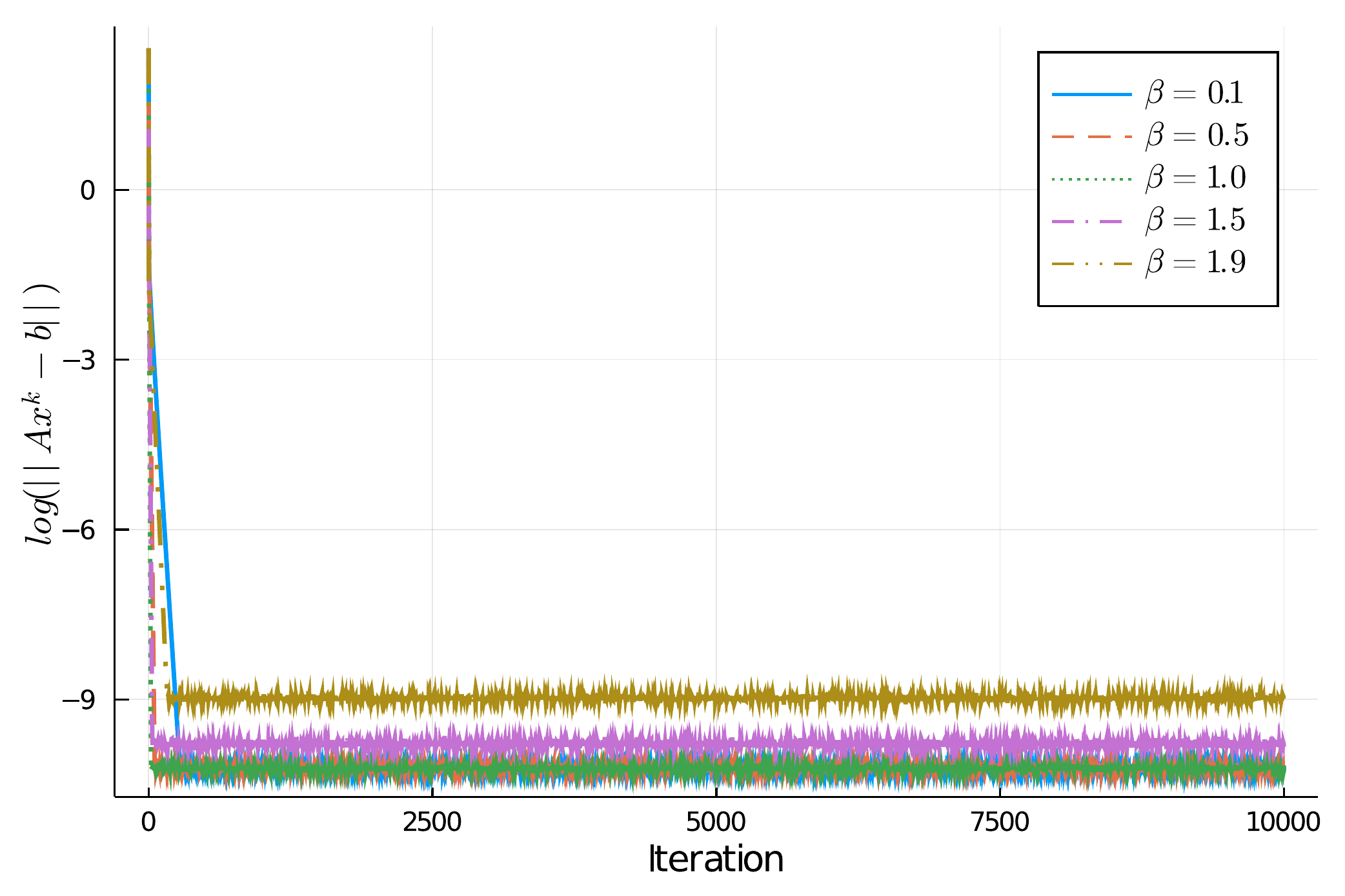} \\[\abovecaptionskip]
    	(a) Infeasibility
  	\end{tabular}
  	\begin{tabular}{@{}c@{}}
    	\includegraphics[width=.4\linewidth]{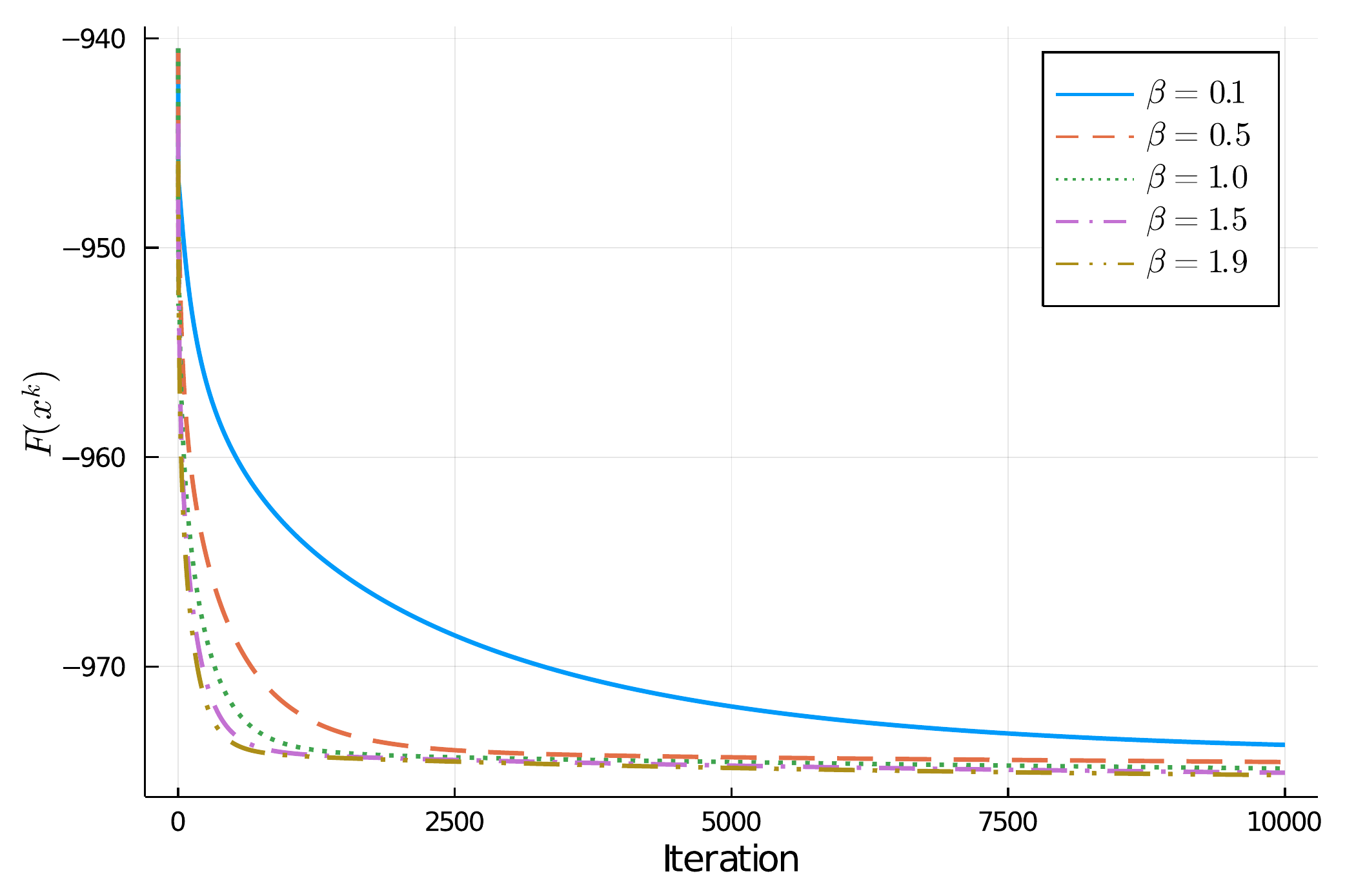} \\[\abovecaptionskip]
    	(b) Objective Value
  	\end{tabular}}
\end{figure*}

Recall the definition of $c_1,c_2$ in \eqref{eq: c_12} and $c_3, c_4$ in \eqref{eq: c_34}. For LCDC-ALM, we set $\nu = \min\{c_1,c_2\}$ and $\rho = 10 \max\{ c_3/(c_1-\nu/2), 2c_3/\nu, 2c_4/\nu\}$. We perform LCDC-ALM with different values of $\beta$, and plot the infeasibility and objective value as functions of iteration in Fig. \ref{fig: compare beta}. The infeasibility converges with similar speed, while a mild slowdown is observed as $\beta$ moves towards its upper or lower bound. In contrast, the objective converges faster when a larger $\beta$ is used. 

\begin{figure*}[h!]
\caption{Comparison of LCDC-ALM (Algorithm \ref{alg: smooth lcdc-alm}) and Proximal ALM \citep{zhang2020global}}
\center{
\label{fig: compare proximal ALM}
	\begin{tabular}{@{}c@{}}
    	\includegraphics[width=.4\linewidth]{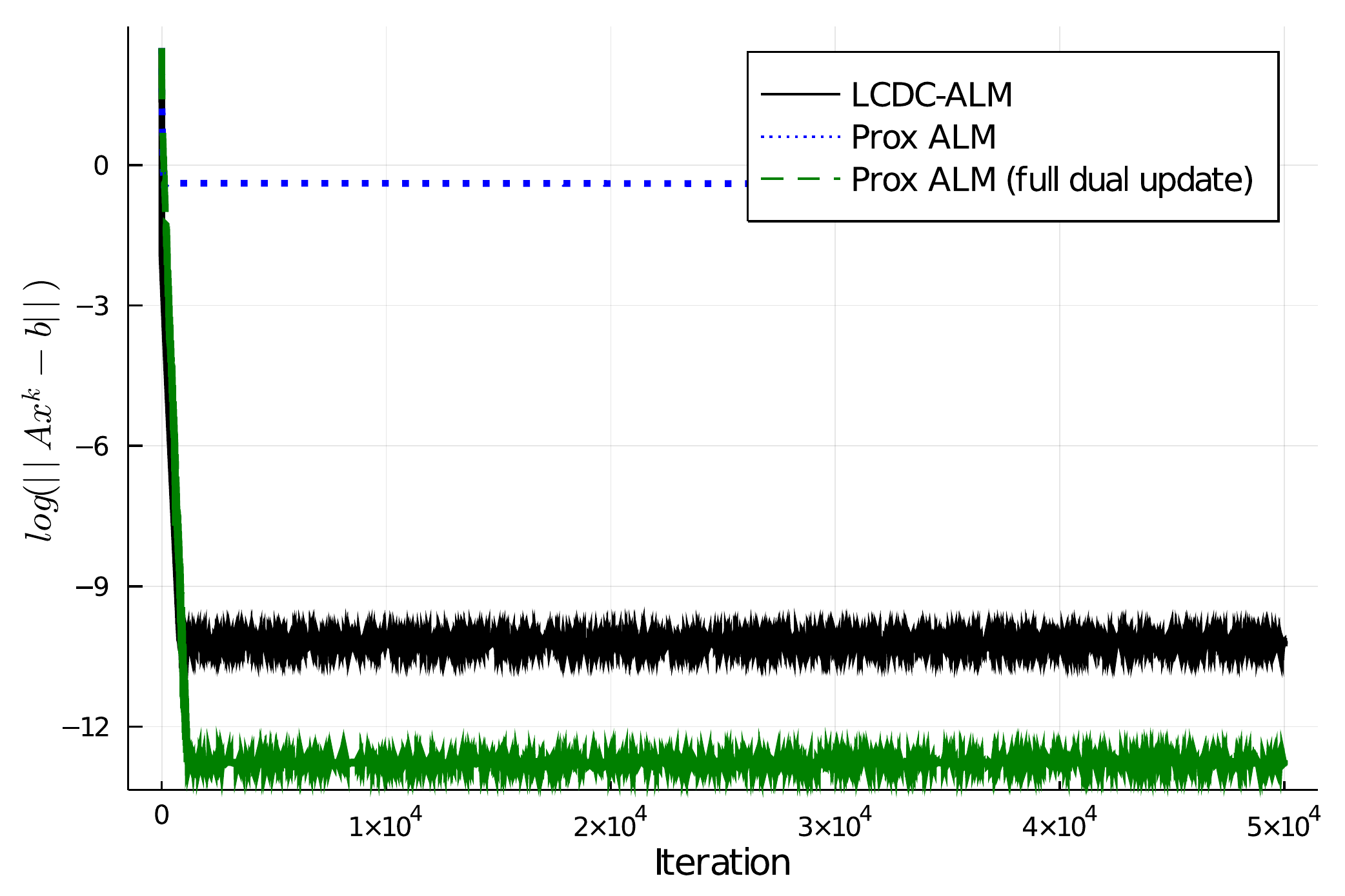} \\[\abovecaptionskip]
    	(a) Infeasibility
  	\end{tabular}
  	\begin{tabular}{@{}c@{}}
    	\includegraphics[width=.4\linewidth]{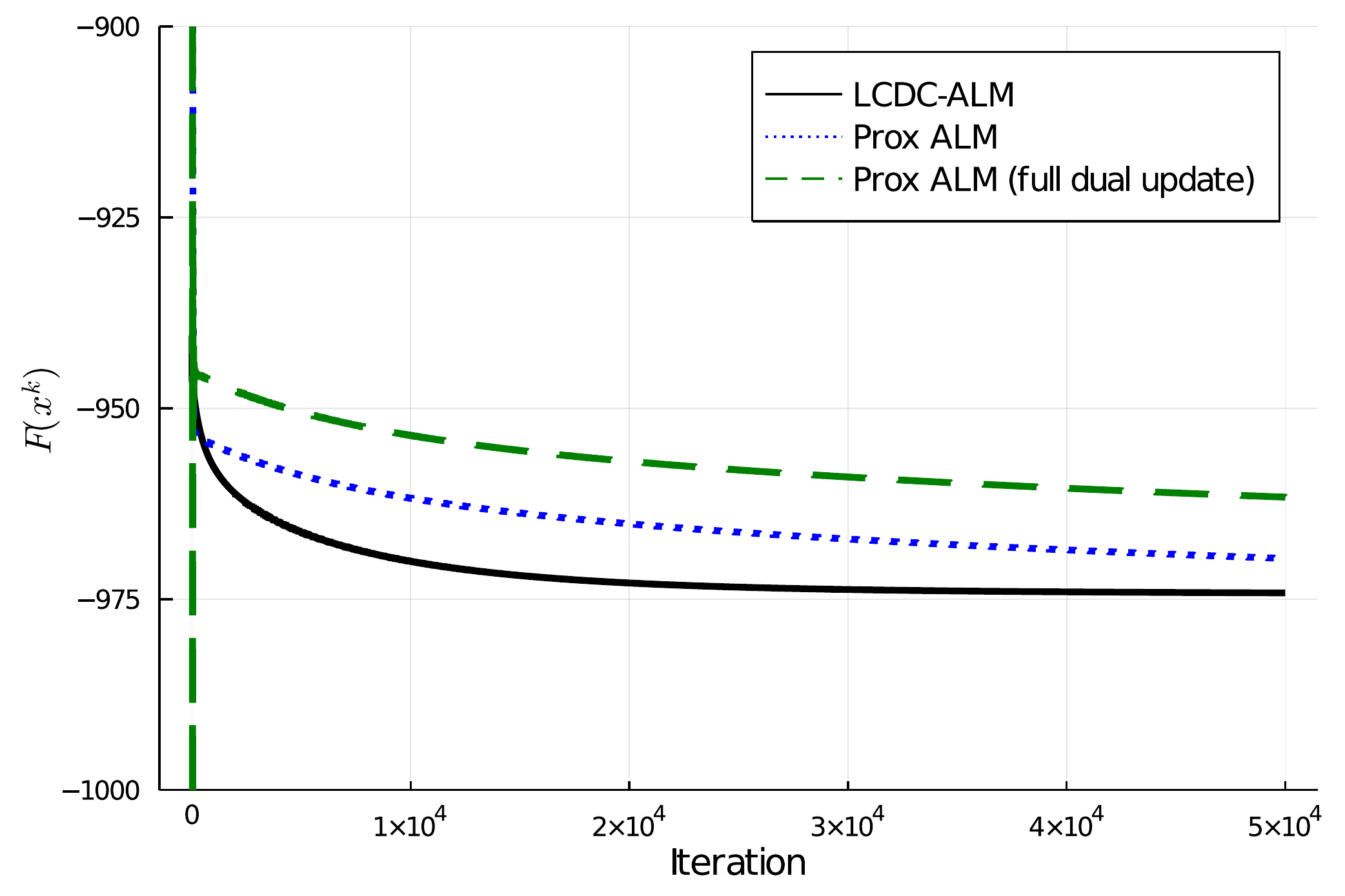} \\[\abovecaptionskip]
    	(b) Objective Value
  	\end{tabular}
  	}
\end{figure*}

Next we compare LCDC-ALM with the proximal ALM proposed in \cite{zhang2020global} in Figure \ref{fig: compare proximal ALM}. For proximal ALM, all parameters are chosen according to \cite[Lemma 3.1]{zhang2020global}: in particular, we choose $\beta = 1/30$ (also for LCDC-ALM), and $\alpha = L_f/[(\|A^\top A\|+4)\|A^\top A\|]$, where $\alpha$ is the dual step size, i.e., $\lambda^{k+1} = \lambda^k + \alpha (Ax^{k+1}-b)$. For the generated instance, we realize the dual step size $\alpha$ is extremely small, and this causes proximal ALM to converge very slowly in both infeasibility and objective. Therefore we also apply a full dual update as $\lambda^{k+1} = \lambda^k+\rho (Ax^{k+1}-b)$, where $\rho = \|Q-G\| =  49.74$. The infeasibility level of the proximal ALM with full dual update then decreases even slightly faster than LCDC-ALM; however, such behavior of proximal ALM is not explained by the analysis in \cite{zhang2020global} and deserves further investigation. Meanwhile, LCDC-ALM achieves a better objective value than proximal ALM.
\section{Concluding Remarks}\label{sec: conclusion}
In this paper, we study the minimization of a DC function $F$ in the form of \eqref{eq: dc} or \eqref{eq: lcdc}. Our algorithmic developments are based on the {the difference-of-Moreau-envelopes} smooth approximation $F_\mu$ introduced in \eqref{eq: smooth approx}. We first study some important properties of $F_{\mu}$, such as Lipschitz smoothness, the correspondence of stationary points, local, and global minima with $F$, and the {coercivity and} level-boundedness of $F_\mu$. Then, we propose algorithms based on the smoothed objective function $F_{\mu}$. In particular, {we show that applying gradient-based updates on $F_\mu$} converges to a stationary solution of $F$ with rate $\mathcal{O}(K^{-1/2})$. Since a local or global solution of $F_{\mu}$ can be used to construct a counterpart of $F$, future directions include exploiting high-order information or sharpness/error bound conditions of $F_\mu$. 

When the minimization of $F$ is explicitly constrained in an affine subspace, we apply the smoothing technique to the classic augmented Lagrangian function and propose two ALM-based algorithms, LCDC-ALM and composite LCDC-ALM, that will find an $\epsilon$-stationary solution in $\mathcal{O}(\epsilon^{-2})$ iterations. We note that due to the more challenging DC setting, the subproblem oracle in composite LCDC-ALM in general cannot be replaced by a single projected gradient step as in \cite{zhang2020global,zhang2020proximal}. We are interested in simplifying ALM subproblems {and extending the smoothing idea to handle linear inequality or even DC constraints} in future works. 

\bibliographystyle{informs2014}
\bibliography{ref}

\newpage

\ECSwitch

\ECHead{Electronic Companion}
\section{Preliminaries of Difference-of-Moreau-Envelopes Smoothing} \label{section: preliminaries}
\subsection{Moreau Envelope}
We summarize some known properties of Moreau envelope in the next proposition. 
\begin{proposition}\label{prop: moreau property}
	Suppose Assumption \ref{assumption: dc} holds and $0 < \mu < 1/m_{\phi}$ in \eqref{eq: phi moreau}. Then the following claims hold.
	\begin{enumerate}
		\item $x_{\mu \phi}$ is Lipschitz continuous with modulus $\frac{1}{1-\mu m_{\phi}}$.
		\item $M_{\mu {\phi}}$ is differentiable with gradient $\nabla M_{\mu \phi}(z) = \mu^{-1}(z-x_{\mu \phi}(z)).$
		\item $\nabla M_{\mu \phi}$ is Lipschitz continuous with modulus $\frac{2-\mu m_{\phi}}{\mu-\mu^2m_{\phi}}$.
	\end{enumerate}
\end{proposition}	
\proof{Proof.}
	The first two claims are well-known, see, e.g., \cite[Lemma 3.5]{zhang2020proximal} and \cite[Proposition 13.37]{rockafellar2009variational}, combining which proves the last one. \Halmos
\endproof

Proposition \ref{prop: moreau property} suggests that $M_{\mu \phi}$ forms a smooth approximation of the possibly nonconvex nonsmooth function $\phi$. Similarly, the Moreau envelope and proximal mapping of $g$ are given by $M_{\mu g}$ and $x_{\mu g}$, respectively. Since $g$ is convex, it is known that $x_{\mu g}$ is 1-Lipschitz and $M_{\mu g}$ is differentiable, whose gradient $\nabla M_{\mu g}(z) = (z-x_{\mu g}(z))/\mu$ is $1/\mu$-Lipschitz \cite[Theorem 6.60]{beck2017first}.

\subsection{Lipschitz Differentiability of $F_{\mu}$}
\begin{proposition}\label{prop: F_mu properties}
Suppose Assumption \ref{assumption: dc} holds and $0 < \mu < 1/m_{\phi}$ in \eqref{eq: phi moreau}. $F_\mu$ is differentiable, and $\nabla F_{\mu}(z) = \mu^{-1}(x_{\mu g}(z)-x_{\mu \phi}(z))$ is Lipschitz continuous with modulus $L_{F_{\mu}}=\frac{2-\mu m_{\phi}}{\mu-\mu^2m_{\phi}}$.
\end{proposition}
\proof{Proof.}
By the Lipschitz differentiability of the Moreau envelope shown in Proposition \ref{prop: moreau property} and the definition of $F_\mu$, we know that $F_{\mu}$ is differentiable, and 
	$\nabla F_\mu(z) = \nabla M_{\mu \phi}(z) - \nabla M_{\mu g}(z) = \mu^{-1}(x_{\mu g}(z)-x_{\mu \phi}(z)).$
	Since $x_{\mu g}$ and $x_{\mu \phi}$ are Lipschitz with modulus 1 and $\frac{1}{1-\mu m_{\phi}}$, respectively, we obtain the claimed $L_{F_{\mu}} = \frac{1}{\mu}(1+\frac{1}{1-\mu m_{\phi}})$. \Halmos
\endproof

If $\phi$ is convex, then the Lipschitz constant of $\nabla F_\mu$ can be improved to $2/\mu$  \citep{hiriart1991regularize}.

\subsection{Correspondence of Stationary Points and Global Minima of $F$ and $F_{\mu}$}\label{sec: correspondence}
In addition to being smooth, the approximation $F_{\mu}$ captures both the local and global structure of the original function $F$. 
In particular, some properties of $F_\mu$ established in \cite{hiriart1991regularize} are summarized in the next proposition. 
\begin{proposition}[\cite{hiriart1991regularize}]\label{prop: known}
	 Suppose {Assumptions \ref{assumption: dc} and \ref{assumption: F finite optimal} hold} and $0 < \mu < 1/m_{\phi}$ in \eqref{eq: phi moreau}. Then the following claims hold.
	 \begin{enumerate}
	 	\item The set of global minimizers of $F_\mu$, $\argmin F_\mu$, is nonempty, and $F^*=$ $\min_{z\in \R^n} F_{\mu}(z)$.
	 	\item (Correspondence of Stationary Point) If ${z}$ is a stationary point of $F_\mu$, i.e., $\nabla F_{\mu}(z)=0$, then $x_{\mu \phi}(z) = x_{\mu g}(z)$, and $x_{\mu \phi}(z)$ is a stationary point of $F$ in the sense of \eqref{eq: dc stationary}  with $F(x_{\mu \phi}(z))=F_{\mu}(z)$; conversely, if $x\in\R^n$ is a stationary point of $F$, then there exists $z\in \R^n$ such that $\nabla F_{\mu}(z) = 0$, $z = x_\phi(x) = x_g(x)$, and $F_{\mu}(z) = F(x)$.
	 	\item (Correspondence of Global Minima) If ${z} \in \argmin F_{\mu}$, then $x_{\mu\phi}({z}) \in \argmin$ $F$; conversely, if $x \in \argmin F$, then there exists $z\in \argmin F_{\mu}$ such that $x = x_{\mu\phi}(z) =x_{\mu g}(z)$. 	
	 \end{enumerate}
\end{proposition}

\section{Proofs in Section \ref{sec: DC}} 
\subsection{Proof of Proposition \ref{prop: map local} }  \label{proof: map local}
\proof{Proof.}
 By Proposition \ref{prop: known} and Lemma \ref{lemma: bounds on F_mu}, we know $x_{\mu \phi}(\bar{z}) = x_{\mu g}(\bar{z})$, and for all $z$ such that  $\|z-\bar{z}\|\leq r$, we have $F(x_{\mu g}(\bar{z}))=F(x_{\mu \phi}(\bar{z})) = F_\mu(\bar{z}) \leq  F_\mu(z) \leq  F(x_{\mu g}(z))$.
    It suffices to show that, for all $x$ sufficiently close to $x_{\mu g}(\bar{z})$, there exists some $z\in \R^n$ such that $\|z-\bar{z}\|\leq r$ and $x = x_{\mu g}(z)$. In particular, take $\tilde{\nabla}g(x) \in \partial g(x)$ and let $z = x+\mu \tilde{\nabla}g(x)$. By construction, we have $0 = \tilde{\nabla}g(x)+ \mu^{-1}(x-z) \in \partial g(x) +\mu^{-1}(x-z)$, and therefore $x = x_{\mu g}(z)$. It follows that 
   	\begin{align*}
   		\|z-\bar{z}\| 
   		&=  \|x+\mu \tilde{\nabla}g(x) - x_{\mu g}(\bar{z})-\mu \tilde{\nabla}g(x_{\mu g}(\bar{z}))\| \\
   		&\leq  \|x-x_{\mu g}(\bar{z})\|+\mu\|\tilde{\nabla}g(x)-\tilde{\nabla}g(x_{\mu g}(\bar{z}))\| \leq  \begin{cases} 
   					(1+ \mu L_g)\|x-x_{\mu g}(\bar{z})\|\leq r\\
   					r-2\mu M_{\partial g} + 2\mu M_{\partial g} = r,
   		 		\end{cases}
   	\end{align*}
   	where the two cases above correspond to the two claims respectively. \Halmos
\endproof

\subsection{Proof of Proposition \ref{prop: level boundedness}}\label{proof: level boundedness}
\proof{Proof.}
\begin{enumerate}
	\item We prove that $\mathrm{lev}_{\alpha} F_\mu$ is bounded. Without loss of generality, we consider $\alpha \geq F^*$, otherwise $\mathrm{lev}_{\alpha} F_\mu = \emptyset$ by Proposition \ref{prop: known}. Firstly notice that 
	\begin{align}
		 \mathrm{lev}_{\alpha} F_{\mu} \subseteq &  \{z:  M_{\mu \phi}(z) -g(z) \leq \alpha \}
			\subseteq \left\{z: \exists x \ \mathrm{s.t.} \ F(x)+ \frac{1}{2\mu}\|x-z\|^2 - L\|x-z\| \leq  \alpha+M \right\}\label{eq: level set inclusion},
	\end{align}
	where the first inclusion is due to $M_{\mu g}\leq g$, and the second inclusion is due to our assumption on $g$. Using the fact that $L\|x-z\| \leq \frac{tL^2}{2} + \frac{\|x-z\|^2}{2t}$ for any $t>0$ and taking $t = 2\mu$, we have
	\begin{align}\label{eq: lower bound F_xz}
		F(x)+ \frac{1}{2\mu}\|x-z\|^2 - L\|x-z\| \geq F(x)+ \frac{1}{4\mu}\|x-z\|^2 - \mu L^2.
	\end{align}
	Now \eqref{eq: level set inclusion}, \eqref{eq: lower bound F_xz}, and the fact that $F(x)+ \frac{1}{4\mu}\|x-z\|^2 \geq \max\{F(x),F^*+\frac{1}{4\mu}\|x-z\|^2\}$ imply that
	\begin{align}
		 \mathrm{lev}_{\alpha} F_{\mu}  
		  \subseteq &  \left \{z: \|z\|\leq \sqrt{(4\mu)(\alpha +M+\mu L^2-F^*)}+ \max_{x:F(x)\leq \alpha+M+\mu L^2} \|x\|\right \}.\label{eq: level set bounded}
	\end{align}
	Since $\mathrm{lev}_{\alpha+M+\mu L^2} F$ is compact, $\max_{x:F(x)\leq \alpha+M+\mu L^2} \|x\|$ is finite, and hence $\mathrm{lev}_{\alpha} F_{\mu} $ is bounded. 
	\item Suppose $ F= \phi-g \geq \alpha \|\cdot\|+ r$ for some $a\in (0, +\infty)$ and $r\in R$. For a convex function $f$, we use $f^*$ to denote its convex conjugate, i.e., $f^*(z) = \sup_{x} \langle z,x \rangle -f(x)$. By definition, we have
	\begin{align}
		\mu F_\mu = \mu (M_{\mu \phi} - M_{\mu g} \notag) = & \left(\mu g + \frac{1}{2}\|\cdot\|^2\right)^* -\left(\mu \phi + \frac{1}{2}\|\cdot\|^2\right)^*\notag \\
		\geq & \left(\mu g + \frac{1}{2}\|\cdot\|^2\right)^* -\left(\mu g + \mu \alpha \|\cdot\| + \mu r  + \frac{1}{2}\|\cdot\|^2\right)^*,\label{eq: conjugate bound 1}
	\end{align}
	where the inequality uses the fact that $f_1\geq f_2$ implies $f_1^*\leq f_2^*$ for any functions $f_1$ and $f_2$. 
	
	We first consider some properties of the first term in \eqref{eq: conjugate bound 1}. For simplicity, denote $p =  \left(\mu g + \frac{1}{2}\|\cdot\|^2\right)^*$. Since $\mu g + \frac{1}{2}\|\cdot\|^2$ is strongly convex with modulus 1, its conjugate $p$ is convex and has Lipschitz gradient with modulus 1 \cite[Theorem 5.26]{beck2017first}. Moreover, we claim that $p$ is coercive: notice that for any $\bar{\alpha} \in (0, +\infty)$, we have
	\begin{align*}
		p(z) =& \max_{x}  \left  \{ \langle z, x\rangle - \mu g(x)-\frac{1}{2}\|x\|^2 \right \}
		\geq \max_{x:\|x\|\leq \bar{\alpha}} \langle z, x\rangle -\max_{x:\|x\|\leq \alpha}  \left  \{\mu g(x)+\frac{1}{2}\|x\|^2 \right \}=  \bar{\alpha} \|z\| + \bar{r},
	\end{align*}
	where $\bar{r} = -\max_{x:\|x\|\leq \alpha}  \left  \{\mu g(x)+\frac{1}{2}\|x\|^2 \right \}$ is finite, since $\mu g+\frac{1}{2}\|\cdot\|^2$ achieves a finite maximum over the compact set $\{x:\|x\|\leq \alpha\}$.
	
	Next we rewrite the second term in \eqref{eq: conjugate bound 1}: for  any $z\in \R^n$,
	\begin{align}
		  \left(\mu g + \mu \alpha \|\cdot\| + \mu r  + \frac{1}{2}\|\cdot\|^2\right)^*(z)= &  \min_x \Big \{ p(x)+ (\alpha \mu \|\cdot\|+\mu r)^*(z-x) \Big\}\notag \\
		  = & \min_{w:\|w\| \leq \alpha \mu}   p(z-w)  - \mu r,\label{eq: conjugare calculus}
	\end{align}
	where the first equality is due to \cite[{Theorem 4.17}]{beck2017first}, and the second equality uses the following facts: $\|\cdot\|^*	= \delta_{\{x:\|x\|\leq 1\}}$ \cite[Section 4.4.2]{beck2017first}, and $(\alpha \mu \|\cdot\| + \mu r)^*(w) = (\alpha \mu)\|\cdot\|^*(\frac{w}{\alpha \mu})  -\mu r$ \cite[{Theorem 4.13, 4.14}]{beck2017first}. Combining \eqref{eq: conjugate bound 1} and \eqref{eq: conjugare calculus}, we have
	\begin{align}
		\mu F_\mu (z) \geq & p(z) -  \min_{w:\|w\| \leq \alpha \mu}   p(z-w)  + \mu r = \max_{w:\|w\| \leq \alpha \mu} p(z)-p(z-w)   + \mu r \notag \\
		\geq & \max_{w:\|w\| \leq \alpha \mu} \langle \nabla p(z), w\rangle -\frac{1}{2}\|w\|^2 + \mu r 
		\geq  \alpha \mu \|\nabla p(z)\| - \frac{1}{2}\alpha^2 \mu^2 + \mu r,\label{eq: conjugate bound 2}
	\end{align}
	where the second inequality is due to the Lipschitz differentiability of $p$, and the last inequality holds with $w = \alpha \mu \frac{\nabla p(z)}{\|\nabla p(z)\|}$ when $\|\nabla p(z)\|>0$, or any $w$ with $\|w\|=\alpha \mu$ when $\|\nabla p(z)\|=0$.
	Notice that \eqref{eq: conjugate bound 2} further suggests that
	\begin{align}
		\liminf_{\|z\|\rightarrow \infty} F_{\mu}(z) \geq & 	\liminf_{\|z\|\rightarrow \infty} \alpha \|\nabla p(z)\| - \frac{1}{2}\alpha^2 \mu + r
		\geq  \alpha\liminf_{\|z\|\rightarrow \infty} \frac{p(z)-p(0)}{\|z\|} - \frac{1}{2}\alpha^2 \mu + r = + \infty, \label{eq: liminf level bounded}
	\end{align}
	where the second inequality is due to the convexity of $p$:  $\|\nabla p(z)\|\|z\| \geq \nabla p(z)^\top z \geq p(z)-p(0)$, and the last equality is due to $p$ being coercive (see an equivalent characterization in \cite[Definition 3.25]{rockafellar2009variational}). Therefore, \eqref{eq: liminf level bounded} implies that $F_\mu$ is level-bounded.
	
	\item Since $\mathrm{dom} ~\phi$ is compact, there exists $R>0$ such that $\mathrm{dom} ~\phi\subseteq \{x:\|x\|\leq R\}$. Notice that 
	\begin{align*}
		 F_\mu(z)  
		\geq & \min_{x:\|x\|\leq R}  \left  \{\phi(x) + \frac{1}{2\mu} \|x\|^2 - \frac{1}{\mu} \langle x, z\rangle \right \} + \max_{x}  \left  \{\frac{1}{\mu} \langle x, z\rangle - g(x)-\frac{1}{2\mu}\|x\|^2 \right \}\\
		\geq &\hat{\phi}^* - \frac{R}{\mu}\|z\| + \max_{x}  \left  \{\frac{1}{\mu} \langle x, z\rangle - g(x)-\frac{1}{2\mu}\|x\|^2\right \},
	\end{align*}
	where in the last inequality, $\hat{\phi}^* = \min_x \{\phi(x) + \frac{1}{2\mu} \|x\|^2 \}$ is well-defined by the strong convexity of $\phi + \frac{1}{2\mu}\|\cdot\|^2$. Pick any $\alpha \in (0, +\infty)$. By a similar argument used in the previous part, we have
	\begin{align*}
		 & \max_{x} \left \{\frac{1}{\mu} \langle x, z\rangle - g(x)-\frac{1}{2\mu}\|x\|^2 \right\} 
		\geq \left(\frac{R}{\mu} + \alpha\right)\|z\|-  \max_{x:\|x\|\leq R+\mu\alpha} \left \{g(x)+\frac{1}{2\mu}\|x\|^2 \right\}.
	\end{align*}
	Combining the above two inequalities, we have $F_\mu(z) \geq \alpha \|z\| + r$, where 
	$r=\hat{\phi}^* - \max_{x:\|x\|\leq R+\mu\alpha} \left \{g(x)+\frac{1}{2\mu}\|x\|^2 \right\}$ is finite. Therefore, we conclude that $F_\mu$ is coercive, and hence also level-bounded.
	\item We show that $\mathrm{lev}_{\alpha} F_{\mu}$ is bounded. Let $z\in \mathrm{lev}_{\alpha} F_{\mu}$. By Lemma \ref{lemma: bounds on F_mu} and the assumption that $F$ is level-bounded, we know that $x_{\mu \phi}(z) \in \mathrm{lev}_{\alpha} F$ and hence is bounded. The definition of $x_{\mu \phi}$ gives $z \in \partial ( \mu \phi + \frac{1}{2}\|\cdot\|^2)(x_{\mu \phi}(z))$. Since $x_{\mu \phi}(z)$ is bounded, and $\mu \phi+\frac{1}{2}\|\cdot\|^2$ is a (strongly) convex function whose domain is $\R^n$, we conclude that $z$ is bounded \cite[Theorem 24.7]{rockafellar1970convex}. \Halmos
\end{enumerate}
\endproof

\section{Proofs in Section \ref{sec: DC algorithm}}
\subsection{Proof of Theorem \ref{thm: gd}}\label{appendix: proof_of_gd}
\proof{Proof.}
	Notice that since $\nabla F_{\mu}$ is $L_{F_{\mu}}$-Lipschitz and $\alpha \leq 1/L_{F_{\mu}}$, we have
	\begin{align*}
			F_{\mu}(z^k) - F_{\mu}(z^{k+1})\geq &  \left(\frac{1}{\alpha}-\frac{L_{F_{\mu}}}{2}\right)\|z^{k+1}-z^k\|^2 \geq  \frac{1}{2\alpha}\|z^{k+1}-z^k\|^2 = \frac{\alpha}{2\mu^2}\|x_{\mu g}(z^k)-x_{\mu \phi}(z^k)\|^2.
	\end{align*}
	Summing the above inequality over $k=0,\cdots,K-1$ for some positive integer $K-1$, we have
	\begin{align}\label{eq: gd telescope}
			\sum_{k=0}^{K-1} \|x_{\mu g}(z^k)-x_{\mu \phi}(z^k)\|^2 \leq \frac{2\mu^2}{\alpha} (F_{\mu}(z^0)-F_{\mu}(z^K)) \leq \frac{2\mu^2}{\alpha} (F_{\mu}(z^0)-F^*).
	\end{align}
	Let $\bar{k}=\argmin_{k=0,\cdots,K-1} \|x_{\mu g}(z^k)-x_{\mu \phi}(z^k)\|^2$, then from \eqref{eq: gd telescope} it holds 
	\begin{align}\label{eq: gd bound}
		\|x_{\mu g}(z^{\bar{k}})-x_{\mu \phi}(z^{\bar{k}})\| \leq \left(\frac{2\mu^2 (F_{\mu}(z^0)-F^*)}{\alpha K}\right)^{1/2}.
	\end{align}
	For any $k\in \Z_+$, due to the optimality of the proximal mapping $x_{\mu \phi}(z^k)$ and $x_{\mu g}(z^k)$, we have
	\begin{align}\label{eq: GD subproblem opt}
		\xi^k &
		=  \mu^{-1}(z^k - x_{\mu \phi}(z^k)) - \mu^{-1}(z^k - x_{\mu g}(z^k)) \in   \partial \phi(x_{\mu \phi}(z^k)) - \partial g(x_{\mu g}(z^k)).
	\end{align}
	In view of \eqref{eq: gd bound}, we have \eqref{eq: gd_opt_coniditon} proved due to the claimed upper bound $K$ in \eqref{eq: gd iter bound}.
	Since $F_{\mu}$ is level-bounded and $\{F_{\mu}(z^k)\}_{k\in \N}$ is monotonically non-increasing, we know the sequence $\{z^k\}_{k\in\N}$ is bounded and therefore has at least one limit point $z^*$. Let $\{z^{k_j}\}_{j\in \N}$ denote the subsequence convergent to $z^*$. Since $x_{\mu \phi}$ and $x_{\mu g}$ are continuous, \eqref{eq: gd telescope} implies $x_{\mu \phi}(z^*) = x_{\mu g}(z^*)$.  Since $g$ is continuous, we have $\lim_{j\rightarrow \infty}g(x_{\mu g}(z^{k_j})) = g(x_{\mu g}(z^*))$; in addition,
	\begin{align*}
		\phi(x_{\mu \phi}(z^*)) \leq & \liminf_{j\rightarrow \infty} \phi(x_{\mu \phi}(z^{k_j}))\leq \limsup_{j\rightarrow \infty} \phi(x_{\mu \phi}(z^{k_j}))\\ \leq & \lim_{j\rightarrow \infty}\Big[\phi(x_{\mu \phi}(z^*))+ \frac{1}{2\mu} \|x_{\mu \phi}(z^*) - z^{k_j}\|^2 - \frac{1}{2\mu} \|x_{\mu \phi}(z^{k_j}) - z^{k_j}\|^2\Big] =\phi(x_{\mu \phi}(z^*)),
	\end{align*}
	where the first inequality is due to the lower-semicontinuity of $\phi$ and the last inequality is due to the optimality of $x_{\mu \phi}(z^{k_j})$ in each Moreau envelope evaluation, and therefore we also have $\lim_{j\rightarrow \infty}\phi(x_{\mu \phi}(z^{k_j})) = \phi(x_{\mu \phi}(z^*))$. Taking limit on \eqref{eq: GD subproblem opt} along the subsequence gives \eqref{eq: dc stationary}.  \Halmos
\endproof

\subsection{Proof of Lemma \ref{lemma: igd descent}} \label{proof: igd descent}
\proof{Proof} 
	We first show that the sequence is bounded from below: 
	$$\mathcal{F}(x^{k}, z^k) \geq \min_x f(x) + h(x) + \frac{1}{2\mu}\|x-z^k\|^2 - M_{\mu g}(z^k) = F_{\mu}(z^k) \geq F^*,$$
	where the last inequality is due to Proposition \ref{prop: F_mu properties}.
	Next we show the descent in $x$:
	\begin{align}
		\mathcal{F}(x^{k+1}, z^k) 
		\leq & f(x^k) + \langle \nabla f(x^k), x^{k+1}-x^{k}\rangle + \frac{L_f}{2}\|x^{k+1}-x^{k}\|^2+h(x^{k+1})+ \frac{1}{2\mu} \|x^{k+1}-z^k\|^2-M_{\mu g}(z^{k}) \notag \\
		\leq & f(x^k) + h(x^k) + \frac{1}{2\mu}\|x^{k}-z^k\|^2-M_{\mu g}(z^k) + \left(\frac{L_f}{2}-\frac{1}{2\mu}\right)\|x^{k+1}-x^{k}\|^2\notag \\
		= & \mathcal{F}(x^{k}, z^k) - 	\left(\frac{\mu^{-1}-L_f}{2}\right)\|x^{k+1}-x^{k}\|^2,\label{eq: igd descent in x}
	\end{align}
	where the first inequality is due to the Lipschitz differentiability of $f$ and the second inequality is due to $x^{k+1}$ being the minimizer of some $\mu^{-1}$-strongly convex function. 
	The descent with respect to $z$ is given as:
	\begin{align}
		& \mathcal{F}(x^{k+1},z^k) - \mathcal{F}(x^{k+1},z^{k+1}) \notag  \\
		 = &  \frac{1}{\mu}\left(\frac{1}{\beta}-\frac{1}{2}\right)\|z^{k+1}-z^k\|^2+ M_{\mu g}(z^{k+1})- M_{\mu g}(z^k)-\left \langle \frac{1}{\mu}(z^k - x_{\mu g}(z^k)), z^{k+1}-z^k\right \rangle \notag \\ 
		 \geq & \frac{1}{\mu}\left(\frac{1}{\beta}-\frac{1}{2}\right)\|z^{k+1}-z^k\|^2,\label{eq: igd descent in z}
	\end{align}
	where we replace $x^{k+1}$ by $x_{\mu g}(z^k)+ \frac{1}{\beta}(z^{k+1}-z^k)$ to get the equality, and the inequality is due to $M_{\mu g}$ being convex and $\nabla M_{\mu g}(z^k) = \mu^{-1}(z^k-x_{\mu g}(z^k))$. Combining \eqref{eq: igd descent in x} and \eqref{eq: igd descent in z} gives \eqref{eq: igd descent}. \Halmos
\endproof

\subsection{Proof of Theorem \ref{thm: igd}}\label{proof: igd}
\proof{Proof.}
	By Lemma \ref{lemma: igd descent}, we have
	$$\mathcal{F}(x^k,z^k) - \mathcal{F}(x^{k+1},z^{k+1}) \geq \min\{c_1, c_2\}\left(\|x^{k+1}-x^k\|^2+(\|z^{k+1}-z^k\|^2 \right),$$
	summing which from $k=0$ to some positive integer $K-1$ gives
	\begin{align}\label{eq: igd telescope}
		\sum_{k=0}^{K-1} \left(\|x^{k+1}-x^k\|^2+(\|z^{k+1}-z^k\|^2 \right) \leq \frac{\mathcal{F}(x^0,z^0) - \mathcal{F}(x^{K},z^{K})}{\min \{c_1,c_2\}} \leq \frac{\mathcal{F}(x^0,z^0) -F^*}{\min \{c_1,c_2\}}.
	\end{align}	
	Let $\displaystyle{\bar{k} = \argmin_{k=0,\cdots, K-1} \|x^{k+1}-x^k\|^2+\|z^{k+1}-z^k\|^2}$, then \eqref{eq: igd telescope} implies 
	\begin{align}\label{eq: igd telescope2}
		\max\{\|x^{\bar{k}+1}-x^{\bar{k}}\|, \|z^{\bar{k}+1}-z^{\bar{k}}\| \}\leq \left(\frac{\mathcal{F}(x^0,z^0) -F^*}{{\min\{c_1,c_2\}}K} \right)^{1/2}.
	\end{align}
	Due to the optimality condition of $x^{k+1}$ and $x_{\mu g}(z^k)$, we have
	\begin{align*}
		\xi^{k+1} 
		\in & \nabla f(x^{k+1})+ {\partial} h(x^{k+1}) - \partial g(x_{\mu g}(z^k)) = \partial \phi(x^{k+1}) - \partial g(x_{\mu g}(z^k)),
	\end{align*}
	{which proves \eqref{eq: igd xi}.} Now in view of \eqref{eq: igd telescope2}, we have
	\begin{align*}
		\max\{\|\xi^{\bar{k}+1}\|, \|x_{\mu g}(z^{\bar{k}}) - x^{\bar{k}+1}\|\}
		\leq &  \left(L_f + \frac{\mu^{-1}+1}{\beta}\right) \left( \frac{\mathcal{F}(x^0,z^0)-F^*}{\min \{c_1,c_2\}K}\right)^{1/2},
	\end{align*}
	{which proves \eqref{eq: igd max} and \eqref{eq: igd iter bound}.}
	
	Next we show that if $F_\mu$ is level bounded, then $\{(x^k,z^k)\}_{k\in \N}$ stays bounded. Since $F_{\mu}$ is continuous and level-bounded,  and $\mathcal{F}(x^0,z^0)\geq \mathcal{F}(x^k,z^k)\geq F_{\mu}(z^k)$, we know {that} $z^k$ {stays} in some compact level set of $F_{\mu}$. Since the mapping $x_{\mu g}$ is continuous, $x_{\mu g}(z^k)$ is also bounded. Consequently, $x^{k+1} = x_{\mu g}(z^k)+\frac{1}{\beta}(z^{k+1}-z^k)$ stays bounded for all $k\in\N$. Therefore, {the sequence $\{(x^k,z^k)\}_{k\in\N}$ has a limit point, denoted as $(x^*,z^*)$}.
	Let $\{(x^{k_j}, z^{k_j})\}_{j\in \N}$ be a subsequence {converging to $(x^*,z^*)$}. {Since $\|x^{k_j}-x^{k_j-1}\|\rightarrow 0$ and $\|z^{k_j}-z^{k_j-1}\|\rightarrow 0$, taking limit on $x^{k+1} = x_{\mu g}(z^k)+\frac{1}{\beta}(z^{k+1}-z^k)$ along the subsequence gives $x^* = x_{\mu g}(z^*)$.}
	{Finally the asymptotical convergence follows a similar argument as in the proof of Theorem \ref{thm: gd}.} \Halmos
\endproof

\section{Proofs in Section \ref{sec: LCDC}}
{
\subsection{Proof of Lemma \ref{lemma: v lowerbounded}} \label{proof: v lowerbounded}
\proof{Proof.}
     We first verify condition \eqref{assumption: lower bounded} under the first two conditions. Since $g$ is Lipschitz, we have $-g(y)\geq -g(x) -L_g\|x-y\| \geq -g(x) -\frac{1}{2}\|x-y\|^2 - \frac{L_g^2}{2}$. For $0 <\mu \leq 1$, it follows that 
		\begin{align}
			v(\mu, \rho)\geq 
			& \inf_{x \in \R^n} \left\{f(x) -g(x) + \frac{\rho}{2}\|Ax-b\|^2 \right\}	-\frac{L_g^2}{2}\label{eq: v_mu_rho_lb}.
		\end{align}
    The first case implies that \eqref{eq: v_mu_rho_lb} is finite for any $\rho \geq 0$. For the second case, notice that we can choose $\rho> 0$ big enough so that $\nabla^2 f + \rho A^\top A \succ 0$. Since $-g$ dominates an affine function, the objective in the right-hand side of \eqref{eq: v_mu_rho_lb} is level-bounded, and hence $v(\mu, \rho) > -\infty$. \\
    Next we verify condition \eqref{assumption: lower bounded} for the third case. Denote $\nabla^2 f = F$ and $\nabla^2 g = G$. The Hessian of the objective in $(x,y)$ in the right-hand side of \eqref{assumption: lower bounded} is positive-definite if $\mu < \lambda_{\max}(G)^{-1}$ and its Schur complement 
    \begin{align*}
        S(\mu, \rho):= & F+\rho A^\top A + \frac{1}{\mu}I_n - \left(-\frac{1}{\mu} I_n\right) \left(\frac{1}{\mu}I_n - G\right)^{-1}\left(-\frac{1}{\mu} I_n\right) \\
        =& F+ \rho A^\top A + \frac{1}{\mu}I_n -\frac{1}{\mu^2} \left[ \mu I_n + \mu^2 G(I_n-\mu G)^{-1}\right] =  F+ \rho A^\top A- G(I_n-\mu G)^{-1}
    \end{align*}
     is positive-definite, where we use the Woodbury matrix identity in the second equality, and $I_n \in \R^{n\times n}$ denotes the identity matrix. Since $F \succ 0$ over the null space of $A$ by assumption, we can choose $\rho>0$ large enough so that $F+ \rho A^\top A \succ 0$. Since $\lambda_{\max}(G(I-\mu G)^{-1}) \rightarrow \lambda_{\max}(G)$ as $\mu \rightarrow 0$, we know $S(\mu, \rho) \succ 0$ if the smallest eigenvalue of $F$ over the null space of $A$ is strictly greater than $\lambda_{\max}(G)$.  This completes the proof. \Halmos
\endproof}

\subsection{Proof of Lemma \ref{lemma: smooth lcdc-alm descent}} \label{proof: smooth lcdc-alm descent}
\proof{Proof.}
	Similar to the derivation in \eqref{eq: igd descent in x}-\eqref{eq: igd descent in z}, the descent of $\psi$ in $x$ and $z$ are given as
	\begin{align*}
		\psi(x^{k}, z^k, \lambda^k) -  \psi(x^{k+1}, z^k, \lambda^k) \geq &\ \left(\frac{\mu^{-1}-L_f}{2}\right)\|x^{k+1}-x^k\|^2,\\
		 \psi(x^{k+1}, z^{k}, \lambda^k)-\psi(x^{k+1}, z^{k+1}, \lambda^k) \geq &\ \frac{1}{\mu}\left(\frac{1}{\beta}-\frac{1}{2}\right)\|z^{k+1}-z^k\|^2.
	\end{align*}
	In addition, 
	\begin{align*}
		 & \psi(x^{k+1}, z^{k+1}, \lambda^k)-\psi(x^{k+1}, z^{k+1}, \lambda^{k+1}) =  \langle \lambda^k -\lambda^{k+1}, Ax^{k+1}-b\rangle = -\frac{1}{\rho}\|\lambda^{k+1}-\lambda^k\|^2. 
	\end{align*}
	Adding the above three expressions completes the proof. \Halmos
\endproof

\subsection{Proof of Lemma \ref{lemma: bound dual by primal}} \label{proof: bound dual by primal}
\proof{Proof.}
	For $k\in \N$, the update of $x^{k+1}$ gives $\nabla f(x^k) + A^\top \lambda^{k+1} +\mu^{-1}(x^{k+1}-z^k) = 0,$
	which implies that for $k\in \N$, $A^\top(\lambda^{k+1}-\lambda^k)= \mu^{-1}(x^k-x^{k+1}) + (\nabla f(x^{k-1})-\nabla f(x^k)) + \mu^{-1}(z^k-z^{k-1}).$
	Since $\lambda^{k+1}-\lambda^k = \rho(Ax^{k+1}-b)$ belongs to the column space of $A$, we have 
	\begin{align*}
		{\sigma_{\min}^+(A)} \|\lambda^{k+1}-\lambda^k\| & \leq \|A^\top(\lambda^{k+1}-\lambda^k)\|\leq  \mu^{-1}\|x^{k+1}-x^k\| + L_f\|x^k-x^{k-1}\| + \mu^{-1}\|z^k-z^{k-1}\|,
	\end{align*}
	{where the first inequality is due to the min-max theorem of eigenvalues of a real symmetric matrix.} 
	Dividing both sides by {$\sigma_{\min}^+(A)$} gives the desired inequality. \Halmos
\endproof

\subsection{Proof of Lemma \ref{lemma: smooth lcdc-alm potential function}} \label{proof: smooth lcdc-alm potential function}
\proof{Proof.}
	\begin{enumerate}
		\item By Lemma  \ref{lemma: bound dual by primal}, squaring both sides gives
		\begin{align*}
			\|\lambda^{k+1}-\lambda^k\|^2 \leq c_3\|x^{k+1}-x^k\|^2 + c_4\|x^k-x^{k-1}\|^2 + c_3\|z^k-z^{k-1}\|^2.
		\end{align*}
		{Then \eqref{eq: potentia descent} follows from Lemma \ref{lemma: smooth lcdc-alm descent} and constants defined in \eqref{eq: kappa}.}
	\item Recall that 
	\begin{align*}
		\Psi_k \geq  {\psi}(x^k, z^k, \lambda^k) 
		\geq & f(x^k)-g(z^k)+ 
		\frac{\rho}{2}\|Ax^k-b\|^2+ \frac{1}{2\mu}\|x^k-z^k\|^2 +  \frac{1}{\rho}\langle \lambda^k, \lambda^{k}-\lambda^{k-1} \rangle \\
		 \geq &  v(\mu, \rho)+ \frac{1}{2\rho}(\|\lambda^k\|^2 - \|\lambda^{k-1}\|^2),
	\end{align*}
	where the second inequality is due to $M_{\mu g}(z) \leq g(z)$ and $\lambda^{k} = \lambda^{k-1}+\rho(Ax^k-b)$, and the third inequality is due to \eqref{assumption: lower bounded}. This further implies that
	\begin{align*}	
			\sum_{k=1}^K (\Psi_k - v(\mu, \rho) )\geq \frac{1}{2\rho}(\|\lambda^K\|^2 - \|\lambda^0\|^2) \geq -\frac{1}{2\rho} \|\lambda^0\|^2 > -\infty,
	\end{align*}
	for all positive integer $K$. Since $\Psi_k$ is non-increasing, we must have $\Psi_k\geq v(\mu, \rho)$ for all $k\in \Z_+$; otherwise, there exists some $\delta>0$ such that $\Psi_k -v(\mu, \rho) < -\delta$ for all large enough $k$, then the above summation would converge to $-\infty$ as $K\rightarrow \infty$. 

	\item For any positive integer $K$, summing \eqref{eq: potentia descent} from 0 to $K-1$ gives 
	\begin{align}
		& \kappa_{\min} K \min_{k=0,\cdots, K-1}\Big\{ \|x^{k+1}-x^k\|^2 +\|z^{k+1}-z^k\|^2  +  \|x^k-x^{k-1}\|^2 + \|z^k-z^{k-1}\|^2\Big\}  \label{eq: telescope} \\
		\leq &\sum_{k=0}^{K-1}\Big(\kappa_1 \|x^{k+1}-x^k\|^2 + \kappa_2\|z^{k+1}-z^k\|^2 + \kappa_3\|x^k-x^{k-1}\|^2 + \kappa_4 \|z^k-z^{k-1}\|^2\Big) \notag \\
		\leq & \sum_{k={0}}^{K-1}\left(\Psi_k-\Psi_{k+1}\right) = \Psi_0 - \Psi_{K} \leq \Psi_0- {v(\mu,\rho)},\notag
	\end{align}
	where the last inequality is due to {$\Psi_K \ge v(\mu,\rho)$} for all $K\ge 1$. Now let $\bar{k}$ be the minimizer in \eqref{eq: telescope}; it follows that
	\begin{align*}
		& \max\left \{\|x^{\bar{k}+1}-x^{\bar{k}}\|^2, \|z^{\bar{k}+1}-z^{\bar{k}}\|^2, {\|x^{\bar{k}}-x^{\bar{k}-1}\|^2}, \|z^{\bar{k}}-z^{\bar{k}-1}\|^2 \right\} \\
		 \leq & \left(\|x^{\bar{k}+1}-x^{\bar{k}}\|^2+  \|z^{\bar{k}+1}-z^{\bar{k}}\|^2+{\|x^{\bar{k}}-x^{\bar{k}-1}\|^2}+\|z^{\bar{k}}-z^{\bar{k}-1}\|^2 \right) \leq  \frac{\Psi_0-{v(\mu,\rho)}}{\kappa_{\min} K}.
	\end{align*}
	This completes the proof. \Halmos
	\end{enumerate}
\endproof
\subsection{Proof of Lemma \ref{lemma: ns_lcdc_descent}} \label{proof: ns_lcdc_descent}
\proof{Proof.}
{By \eqref{eq: composite-lcdc-primal-update} and the $\mu$-strong convexity of the function in the following line, it holds for all $x\in \mathcal{H}$ that
\begin{align}
	& \langle \nabla f(x^k)-\xi^k_g, x-x^k \rangle + h(x) + \langle \lambda^k, Ax-b\rangle + \frac{\rho}{2}\|Ax-b\|^2 + \frac{1}{2\mu} \|x-z^k\|^2 \notag \\
	\geq &  \langle \nabla f(x^k)-\xi^k_g, x^{k+1}-x^k \rangle + h(x^{k+1}) + \langle \lambda^k, Ax^{k+1}-b\rangle + \frac{\rho}{2}\|Ax^{k+1}-b\|^2 \notag  \\
		& +  \frac{1}{2\mu} \|x^{k+1}-z^k\|^2 + \langle \zeta^{k+1}, x-x^{k+1}\rangle + \frac{1}{2\mu}\|x^{k+1}-x\|^2. \label{eq: strong cvx composite}
\end{align}
Using the Lipschitz condition of $\nabla f$, the convexity of $g$, inequality \eqref{eq: strong cvx composite} with $x=x^k$, the fact that $\langle \zeta^{k+1}, x^{k+1}-x^k\rangle\leq \frac{1}{4\mu}\|x^{k+1}-x^k\|^2 + \mu \|\zeta^{k+1}\|^2$, and the choice that $\|\zeta^{k+1}\| \leq \epsilon_{k+1}$, the descent in $x$ can be derived as follows:
\begin{align*}
	 P(x^{k+1}, z^k, \lambda^k)  \leq  & P(x^k, z^k, \lambda^k) + \mu \epsilon^2_{k+1} - \frac{\mu^{-1} - 2L_f}{4}\|x^{k+1}-x^k\|^2. 
\end{align*}
The descent with respect to variable $z$ can be derived as follows:
\begin{align}
     P(x^{k+1}, z^{k}, \lambda^k)-P(x^{k+1}, z^{k+1}, \lambda^k)
     = & \frac{1}{2\mu}\left(\frac{2}{\beta}-1\right)\left\|z^{k+1}-z^k\right\|^2 \geq \frac{1}{2\beta\mu}\left\|z^{k+1}-z^k\right\|^2, \notag 
\end{align}
where we replace $x^{k+1}$ by $z^k + \frac{1}{\beta}({z}^{k+1}-z^k)$ to get the equality, and the inequality is due to $\beta \leq 1$.
Finally, similar to Lemma \ref{lemma: smooth lcdc-alm descent}, the change in $\lambda$ is  
\begin{align}\notag 
     P(x^{k+1}, z^{k+1}, \lambda^k)-P(x^{k+1}, z^{k+1}, \lambda^{k+1})= -\frac{1}{\rho}\|\lambda^{k+1}-\lambda^k\|^2.
\end{align}
Combining the above three expressions completes the proof. \Halmos}
\endproof

\subsection{Proof of Lemma \ref{lemma: bound init val}} \label{proof: bound init val}
\proof{Proof.}
{The optimality condition \eqref{eq: strong cvx composite} with $k=0$ and $x = \bar{x}$ (by Assumption \ref{assumption: nonsmooth lcdc}, $A\bar{x}=b$) gives
\begin{align}
    & \langle \nabla f(x^{0})-\xi_g^{0}, x^{1}-x^{0} \rangle  + h(x^1) +  \langle \lambda^{0}, Ax^{1}-b\rangle + \frac{\rho}{2}\|Ax^1-b\|^2 + \frac{1}{2\mu}\|x^1-z^{0}\|^2 	\notag \\
    \leq & \langle \nabla f(x^{0})-\xi_g^{0}, \bar{x}-x^{0} \rangle  + h(\bar{x}) + \frac{1}{2\mu}\|\bar{x}-z^{0}\|^2 + \langle \zeta^1, x^{1}-\bar{x} \rangle \label{eq: ns_lcdc_alm_opt_x}.
\end{align}
Since $\|\zeta^{1}\|\leq \epsilon_1 \leq 1$, $\mu < L_f^{-1}$, and $x^1, \bar{x}\in \mathcal{H}$, we have
\begin{align}
    \langle \zeta^1, x^{1}-\bar{x} \rangle \leq \frac{\mu}{2}\|\zeta^1\|^2 + \frac{1}{2\mu}\|x^{1}-\bar{x} \|^2\leq \frac{L_f^{-1}}{2} + \frac{1}{2\mu}D_\mathcal{H}^2. \label{eq: bound inner prod}
\end{align}
The above two inequalities together  with the $L_h$-Lipschitz continuity of $h$ and $\bar{x},z^0\in\cal{H}$ give that 
\begin{align}\label{eq: ns_lcdc_alm_bound_initial_infeas}
    \rho \|Ax^1-b\|^2 \leq 2(M_{\nabla f}+M_{\partial g} + L_h)D_{\mathcal{H}} + 2\mu^{-1}D_{\mathcal{H}}^2 +L_f^{-1}+ 2\|\lambda^{0}\|\max_{x\in \mathcal{H}}\|Ax-b\|.
\end{align}
Notice that due to $\nabla f$ being Lipschitz and $g$ being convex, $P(x^1, z^{0}, \lambda^{0})$ is bounded from above by 
	\begin{align}\label{eq: ns_lcdc_alm_bound_initial_potential}
		 & f(x^{0}) -g(x^{0})+ \frac{L_f}{2}\|x^1-x^{0}\|^2 + \langle \nabla f(x^{0})-\xi_g^{0}, x^{1}-x^{0} \rangle  + h(x^1) +  \langle \lambda^{0}, Ax^{1}-b\rangle + \frac{\rho}{2}\|Ax^1-b\|^2 + \frac{1}{2\mu}\|x^1-z^{0}\|^2 \notag \\
		\leq  &  f(x^{0}) -g(x^{0})+ \frac{L_f}{2}\|x^1-x^{0}\|^2 +\langle \nabla f(x^{0})-\xi_g^{0}, \bar{x}-x^{0} \rangle  + h(\bar{x}) + \frac{1}{2\mu}\|\bar{x}-z^{0}\|^2  + \langle \zeta^1, x^{1}-\bar{x}\rangle \notag \\
		\leq  & \max_{x\in \mathcal{H}} \{f(x)+h(x)-g(x)\} + (L_h+ M_{\nabla f} + M_{\partial g})D_{\mathcal{H}} + \frac{L_f+ 2 \mu^{-1}}{2} D_{\mathcal{H}}^2 + \frac{L_f^{-1}}{2},
	\end{align}
	where the first inequality is due to \eqref{eq: ns_lcdc_alm_opt_x}, and the second inequality is due to the compactness of $\mathcal{H}$, the $L_h$-Lipschitz continuity of $h$, and \eqref{eq: bound inner prod}. By Lemma \ref{lemma: ns_lcdc_descent}, \eqref{eq: ns_lcdc_alm_bound_initial_infeas}, and \eqref{eq: ns_lcdc_alm_bound_initial_potential}, we have
	\begin{align*}
		P(x^1,z^1, \lambda^1) = 	P(x^1, z^{1}, \lambda^{0}) + \rho\|Ax^{1}-b\|^2 \leq P(x^1, z^{0}, \lambda^{0}) + \rho\|Ax^{1}-b\|^2 \leq \overline{P}.
	\end{align*} 
	This completes the proof. \Halmos }
\endproof

\subsection{Proof of Lemma \ref{lemma: ns dual bd}}\label{proof: ns dual bd}
\proof{Proof.}
By step 4 in Algorithm \ref{alg: ns-lcdc-alm}, there exists $\xi^{k+1}_h \in \partial h(x^{k+1})$ such that 
\begin{align}\label{eq: opt_ns_x}
	{\zeta^{k+1}} = \nabla f(x^{k})-\xi^k_g + \xi_h^{k+1} + A^\top \lambda^{k+1}+ \frac{1}{\mu}(x^{k+1}-z^k).	
\end{align}
Since ${\lambda^{0}}, b \in \mathrm{Im}(A)$ and $\lambda^{k+1} = \lambda^k + \rho(Ax^{k+1}-b)$, $\lambda^{k+1} \in \mathrm{Im}(A)$ for all {$k\in \N$}. 
Since ${z^0} \in \mathcal{H}$, $x^k \in \mathcal{H}$, and $z^{k+1} = (1-\beta) z^k+ \beta x^{k+1}$, $z^k\in \mathcal{H}$ for all $k\in \N$ as well. Consequently, 
\begin{align}\label{eq: dual bound 1}
	 \|\lambda^{k+1}\|\leq \frac{1}{{\sigma_{\min}^+(A)}}\left(M_{\nabla f}+ M_{\partial g} + \frac{D_{\mathcal{H}}}{\mu} {+ 1} \right) + \frac{ \|\xi_h^{k+1}\|}{{\sigma_{\min}^+(A)}}. 
\end{align}
By \cite[Lemma 4.7]{melo2020iteration}, we can bound $\|\xi_h^{k+1}\|$ as follows:
\begin{align}\label{eq: renato lemma}
	\bar{d} \|\xi_h^{k+1}\| \leq (\bar{d} + \|x^{k+1}-\bar{x}\|)L_h +  \left\langle \xi_h^{k+1}, x^{k+1}-\bar{x}\right\rangle \leq 2D_{\mathcal{H}}L_h + \left\langle \xi_h^{k+1}, x^{k+1}-\bar{x}\right\rangle.
\end{align}
Using \eqref{eq: opt_ns_x}, we can further bound the inner product term in \eqref{eq: renato lemma} by 
\begin{align}\label{eq: dual bound 2}
	 \langle \xi_h^{k+1}, x^{k+1}-\bar{x}\rangle 
	& =  \left\langle {\zeta^{k+1}} -\nabla f(x^{k})+\xi^k_g - A^\top \lambda^{k+1} - \frac{1}{\mu}(x^{k+1}-z^k), 
	 x^{k+1} - \bar{x}\right \rangle \notag \\
	& \leq  \left (M_{\nabla f} + M_{\partial g} + \frac{D_{\mathcal{H}}}{\mu} {+1} \right) D_{\mathcal{H}} - \left\langle \lambda^{k+1}, Ax^{k+1}-b \right\rangle \notag \\ 
	& {\leq}   \left (M_{\nabla f} + M_{\partial g}  + \frac{D_{\mathcal{H}}}{\mu} {+1} \right) D_{\mathcal{H}}+\frac{1}{\rho}\|\lambda^{k+1}\|\|\lambda^k\|-\frac{1}{\rho}\|\lambda^{k+1}\|^2,
\end{align}
where we use the facts that $A\bar{x}=b$ {and $\|\zeta^{k+1}\|\leq \epsilon_{k+1}\leq 1$} in the first inequality, and $Ax^{k+1}-b = \frac{1}{\rho}(\lambda^{k+1}-\lambda^k)$ to get the second {inequality}. Combining \eqref{eq: dual bound 1}, \eqref{eq: renato lemma} and \eqref{eq: dual bound 2}, we have 
\begin{align*}
	\frac{\|\lambda^{k+1}\|^2}{\rho {\sigma_{\min}^+(A)} } + \bar{d}\|\lambda^{k+1}\|
	 \leq & \frac{\|\lambda^{k+1}\| \|\lambda^k\|}{\rho {\sigma_{\min}^+(A)} }+\frac{2D_\mathcal{H}}{{\sigma_{\min}^+(A)} }\left(M_{\nabla f}+ M_{\partial g} +\frac{D_{\mathcal{H}}}{\mu}+L_h {+1} \right),
\end{align*}
which further implies that, {for all $k\in \N$},
\begin{align*}
	\left(\frac{\|\lambda^{k+1}\|}{\rho {\sigma_{\min}^+(A)} } + \bar{d} \right)\|\lambda^{k+1}\| \leq \frac{\|\lambda^{k+1}\|}{\rho {\sigma_{\min}^+(A)} } \|\lambda^k\| + \bar{d} \Lambda,
\end{align*}
The claim then follows from an inductive argument: {if $\|\lambda^k\| \leq \Lambda$, then the above inequality implies that $\|\lambda^{k+1}\| \leq \Lambda$ as well, and $\|\lambda^0\|\leq \Lambda$ holds by the definition of $\Lambda$. }  \Halmos
\endproof

\subsection{Proof of Lemma \ref{lemma: ns_lcdc}} \label{proof: ns_lcdc}
\proof{Proof.}
\begin{enumerate}
	\item For all $k\in \N$, 
	\begin{align*}
		 P(x^k, z^k, \lambda^k) = & f(x^k)+h(x^k)-g(x^k) + \langle \lambda, Ax^k-b\rangle + \frac{\rho}{2}\|Ax^k-b\|^2 + \frac{1}{2\mu}\|x^k-z^k\|^2\\
		\geq &  \min_{x\in \R^n} \left\{f(x)+h(x)-g(x)\right\} - \Lambda \max_{x\in \mathcal{H}}\|Ax-b\| > -\infty,
	\end{align*}
	where the inequality is due to the continuity of $f$, $h$, $g$, and $\|\cdot\|$ over compact domain $\mathcal{H}$.
	\item Lemma \ref{lemma: ns_lcdc_descent} implies that, for all $k\in \N$, 
	\begin{align*}
		  \eta \left( \|x^{k+1}-x^k\|^2 + \|z^{k+1}-z^k\|^2\right)
		\leq & P(x^k,z^k,\lambda^k) - P(x^{k+1},z^{k+1},\lambda^{k+1}) + \frac{1}{\rho}\|\lambda^{k+1}-\lambda^k\|^2 {+ \mu \epsilon_{k+1}^2}.
	\end{align*}
	Summing the above inequality over $k={1},\dots, {K}$,
	\begin{align*}
		& \eta K \min_{k\in [K]} \left\{ \|x^{k+1}-x^k\|^2 + \|z^{k+1}-z^k\|^2\right\} \leq   \eta \sum_{k=1}^{K}\left\{\|x^{k+1}-x^k\|^2 + \|z^{k+1}-z^k\|^2\right\}\\
		\leq & P(x^1, z^1, \lambda^1) -  P(x^{K+1}, z^{K+1}, \lambda^{K+1})+ \frac{1}{\rho} \sum_{k=1}^{K}\|\lambda^{k+1}-\lambda^k\|^2 + {\mu \sum_{k=1}^K \epsilon_{k+1}^2 } \\
		\leq & {\overline{P}}-\underline{P} {+\mu E} + \frac{4}{\rho} \sum_{k=1}^{K+1} \|\lambda^k\|^2 \leq {\overline{P}}-\underline{P}  {+\mu E} +   \frac{4(K+1)\Lambda^2}{\rho}, 
	\end{align*}
	where we use the the fact that $P(x^{{K+1}}, z^{{K+1}}, \lambda^{{K+1}})\geq \underline{P}$ {for all integer $K\in \Z_+$} and $\|\lambda^{k+1}-\lambda^k\|^2\leq 2\|\lambda^{k+1}\|^2 + 2\|\lambda^k\|^2$ in the third inequality, and Lemma \ref{lemma: ns dual bd} in the last inequality. Let $\bar{k}$ be the minimizer in the first line of the above chain of inequalities, then dividing both sides by any positive integer $K$ gives
	\begin{align*}
			  \max \left\{ \|x^{\bar{k}+1} -x^{\bar{k}}\|^2, \|z^{\bar{k}+1} -z^{\bar{k}}\|^2\right\} 
		\leq \frac{{\overline{P}}-\underline{P}  {+\mu E} }{\eta K} + \frac{8 \Lambda^2}{\eta\rho}.
	\end{align*}
	This completes the proof. \Halmos
\end{enumerate}	
\endproof

\end{document}